\documentstyle[11pt,amsmath,amssymb]{article}

\begin{document}
\newtheorem{lem}{Lemma}[section]
\newtheorem{th}{Theorem}[section]
\newtheorem{prop}{Proposition}[section]
\newtheorem{rem}{Remark}[section]
\newtheorem{define}{Definition}[section]
\newtheorem{cor}{Corollary}[section]
\newtheorem{exa}{Example}[section]

\makeatletter\@addtoreset{equation}{section}\makeatother
\def\theequation{\arabic{section}.\arabic{equation}}

\renewcommand{\kappa}{\varkappa}
\newcommand{\ext}{{mk}}
\newcommand{\N}{{\Bbb N}}
\newcommand{\C}{{\Bbb C}}
\newcommand{\Z}{{\Bbb Z}}
\newcommand{\R}{{\Bbb R}}
\newcommand{\Rp}{{\R_+}}
\newcommand{\eps}{\varepsilon}
\newcommand{\Om}{\Omega}
\newcommand{\OM}{\Om_X^{M}}
\newcommand{\om}{\omega}
\newcommand{\PR}{\operatorname{Pr}}
\newcommand{\XM}{{X\times M}}
\newcommand{\supp}{\operatorname{supp}}
\newcommand{\la}{\langle}
\newcommand{\ra}{\rangle}
\newcommand{\Oc}{{\cal O}_{\mathrm{c}}}
\newcommand{\Bc}{{\cal B}_{\mathrm{c}}}
\newcommand{\Lext}{\Lambda_{{\mathrm \ext}}}
\newcommand{\tsigma}{{\tilde \sigma}}
\renewcommand{\emptyset}{\varnothing}
\renewcommand{\tilde}{\widetilde}
\newcommand{\Got}{{\frak G}}
\newcommand{\got}{{\frak g}}
\newcommand{\Diff}{\operatorname{Diff}_0(X)}
\newcommand{\RtX}{\R_+^X}
\newcommand{\rom}[1]{{\rm #1}}
\newcommand{\FC}{{\cal F}C_{\mathrm b}^\infty\big({\frak D},\OM\big)}
\newcommand{\FCO}{{\cal F}C_{\mathrm b}^\infty\big({\frak D}_0,
\OM\big)}
\newcommand{\XR}{X\times\R_+}
\newcommand{\vph}{\varphi}
\newcommand{\fii}{\vph}
\newcommand{\di}{\partial}
\renewcommand{\div}{\operatorname{div}}
\newcommand{\Bhat}{\widehat{{\cal B}}}
\newcommand{\un}{{\underline n}}
\newcommand{\pit}{{\pi_\tsigma}}
\newcommand{\FP}{{\cal F}{\cal P}({\frak D},\OM)}
\newcommand{\FCp}{{\cal F}C^\infty_{\mathrm p}({\frak D},\OM)}
\newcommand{\FPO}{{\cal F}{\cal P}({\frak D}_0,\OM)}
\newcommand{\FCpO}{{\cal F}C^\infty_{\mathrm p}({\frak D}_0,\OM)}
\newcommand{\EOM}{{\cal E}^\Omega_{\pi_\tsigma}}
\newcommand{\HOM}{H^\Om_{\pi_\tsigma}}
\newcommand{\HXM}{H_\tsigma^\XM}
\newcommand{\Agot}{{\frak A}}
\newcommand{\nablat}{{\widetilde\nabla}}

\renewcommand{\title}[1]{\leftline{\Large\bf #1}\par\medskip}
\renewcommand{\author}[1]{\medskip{\large #1}\par\medskip}
\newcommand{\address}[1]{\noindent{\sl #1}}

\title{ANALYSIS AND GEOMETRY}
\title{ON MARKED CONFIGURATION SPACES}

\author{SERGIO ALBEVERIO, YURI KONDRATIEV}
\author{EUGENE LYTVYNOV, AND GEORGI US}

\begin{center}\bf Abstract\end{center}

\noindent\begin{small} We carry out analysis and geometry on a
marked configuration space $\Omega^M_X$ over a Riemannian manifold
$X$ with marks from a space $M$. We suppose that $M$ is a
  homogeneous   space $M$ of a Lie
group $G$. As a transformation group $\frak A$ on $\Omega_X^M$
  we take the ``lifting''
to $\Omega_X^M$ of the action on $X\times M$ of the semidirect
product of the group $\operatorname{Diff}_0(X)$ of diffeomorphisms
on $X$ with compact support and the group $G^X$ of smooth
currents, i.e., all $C^\infty$ mappings of $X$ into $G$ which are
equal to the identity element outside of a compact set. The marked
Poisson measure $\pi_\sigma$ on $\Omega_X^M$ with L\'evy measure
$\sigma$ on $X\times M$ is proven to be quasiinvariant under the
action of $\frak A$. Then, we derive a geometry on $\Omega_X^M$ by
a natural ``lifting'' of the corresponding geometry on $X\times
M$. In particular, we construct a gradient $\nabla^\Omega$ and a
divergence $\operatorname{div}^\Omega$. The associated volume
elements, i.e., all probability measures $\mu$ on $\Omega_X^M$
with respect to which $\nabla^\Omega$ and
$\operatorname{div}^\Omega$ become dual operators on
$L^2(\Omega_X^M;\mu)$, are identified as the mixed marked Poisson
measures with mean measure equal to a multiple of $\sigma$. As a
direct consequence of our results, we obtain marked Poisson space
representations of the group $\frak A$ and its Lie algebra $\frak
a$.   We investigate also Dirichlet forms and Dirichlet operators
connected with (mixed) marked Poisson measures. \vspace{5mm}

\setcounter{section}{-1}

\end{small}

\section{Introduction}
In recent  years, stochastic analysis and differential geometry on
configuration spaces have been considerably developed in a series
of papers  [5--8], see also \cite{Roe, ADL1, ADL2}. It has been
shown, in particular,  that the geometry of the configuration
space $\Gamma_X$ over a Riemannian manifold $X$ can be constructed
via a simple ``lifting procedure'' and is completely determined by
the Riemannian structure of $X$. The mixed Poisson measures are
then exhibited as the  ``volume elements'' corresponding to the
differential geometry introduced on $\Gamma_X$. Intrinsic
Dirichlet forms and operators, their canonical processes, as well
as  Gibbs measures on configuration spaces, their characterization
by  integration by parts, and the corresponding  stochastic
dynamics  are among the problems which have been treated in the
above framework.

A starting point for this analysis, more exactly, for the
definition of differentiation on the configuration space, was the
representation of the group of diffeomorphisms $\Diff$ on $X$ with
compact support that was constructed by G.~A.~Goldin et~al.\
\cite{zabyl} and A.~M.~Vershik et~al.\ \cite{VGG2} (see also
\cite{Ob,Sh,Ismagilov}). The construction of this representation
used, in turn, the fact, following from the Skorokhod theorem,
that the Poisson measure is quasiinvariant with respect to the
group $\Diff$.

On the other hand, starting with the same work \cite{VGG2}, many
researchers consider representations also on marked (in
particular, compound) Poisson spaces. In statistical mechanics of
continuous systems, marked Poisson measures and their Gibbsian
perturbations are used for the description of many concrete
models, see e.g.\ \cite{AGL78}. Hence, it is natural to ask about
geometry and analysis on marked Poisson spaces. The first work in
this direction was the paper \cite{KSS}, in which, just as in the
case of the usual Poisson measure, the action of the group $\Diff$
was used for the definition of the differentiation. However, this
group proved to be too small for reconstructing mixed marked
Poisson measures as ``volume elements,'' which means that $\Diff$
is to be extended in a proper way, which we will describe in the
present paper.

Let us recall that the configuration space $\Gamma_X$ is defined
as the space of all locally finite subsets (configurations) in
$X$. Then, the marked configuration space $\Omega_X^M$ over $X$
with marks from, generally speaking, a manifold $M$ is defined as
\[ \Omega_X^M:=\big\{\,(\gamma,s)\mid \gamma\in \Gamma_X,\, s\in M^\gamma\,\big\},\]
where     $M^\gamma$ stands for the set of all maps $\gamma\ni
x\mapsto s_x\in M$. Let $\tsigma$ be a Radon measure on $ X\times
M$ such that $\tsigma(K\times M)<\infty$ for each compact
$K\subset X$ and $\tsigma$ is nonatomic in $X$, i.e.,
$\tsigma(\{x\}\times M)=0$ for each $x\in X$. Then, one can define
on $\Omega_X^M$ a marked Poisson measure $\pi_\tsigma$ with L\'evy
measure $\tsigma$.

Of course, one could consider $\pi_\tsigma$ as a usual Poisson
measure on  the configuration space $\Gamma_{X\times M}$ over the
Cartesian product of the underlying manifold $X$ and the space of
marks $M$, and study the properties of this measure using the
results of [2--5]. However, such an approach does not distinguish
between the two different natures of $X$ and $M$ and the different
roles that these play in physics. Thus, our aim is to introduce
and study such transformations of the marked configuration space
which do ``feel'' this difference and lead to an appropriate
stochastic analysis and differential geometry.

In our previous paper \cite{KLU1}, we were concerned with the
model case  $M=\R_+$, which corresponds, in fact, to the case of a
compound Poisson measure. As has been promised in \cite{KLU1}, we
generalize in the present paper the results of \cite{KLU1} to the
case where $M$ is a homogeneous space of a Lie group $G$. This
situation is natural from the physical point of view. For example,
one can take $X=\R^3$ and $M$ to be the unit sphere $S^2$ in
$\R^3$, and consider any marked configuration
$(\gamma,s)=\{(x,s_x)_{x\in\gamma}\}\in\Omega_X^M$ as a system of
particles in $\R^3$ situated at the points $x$  of $\gamma$ and
having spin $s_x$ at $x\in\gamma$. One has then to take $G$ as the
rotation group, see e.g.\ \cite{GMS}.

Let $G^X$ denote the group of smooth currents, i.e., all
$C^\infty$ mappings $X\ni x\mapsto \eta(x)\in G$ which are equal
to the identity element of $G$ outside of a compact set (depending
on $\eta$). We define the group $\frak A$  as the semidirect
product of the groups $\Diff$ and $G^X$: for $a_1=(\psi_1,\eta_1)$
and $a_2=(\psi_2, \eta_2)$, where $\psi_1,\psi_2\in \Diff$ and
$\eta_1,\eta_2\in G^X$, the multiplication of $a_1$ and $a_2$ is
given by
\[ a_1a_2=(\psi_1\circ\psi_2,\eta_1(\eta_2\circ\psi_1^{-1})).\]
The group $\frak A$ acts in $\XM$ as follows: for any $a=(\psi,
\eta)\in{\frak A}$
\[ \XM\ni(x,m)\mapsto a(x,m)=(\psi(x),\eta(\psi(x))m)\in\XM,\]
where,  for $g\in G$ and $m\in M$, $gm$ denotes the action of $g$
on $m$. Since each $\om\in\OM$ can be interpreted as a subset of
$\XM$, the action of $\frak A$ can be lifted to an action in
$\OM$. The marked Poisson measure $\pi_\tsigma$ is proven to be
quasiinvariant under it. Thus, we can easily construct, in
particular, a representation of $\frak A$ in $L^2(\pi_\tsigma)$.
It should be stressed, however, that our representation of $\frak
A $ is reducible, because so is the regular representation of
$\frak A$ in $L^2(\tsigma)$, see subsec.~3.5 in \cite{KLU1} for
details.

Having introduced the action of the group $\frak A$ on $\OM$, we
proceed to derive analysis and geometry on $\OM$ in a way parallel
to the works \cite{AKR, KLU1}, dealing with the usual
configuration space $\Gamma_X$ and the marked configuration space
$\Omega^{\R_+}_X$, respectively. In particular, we note that the
Lie algebra $\frak a$ of the group $\frak A$ is given by ${\frak
a} =V_0(X)\times C_0^\infty(X;{\frak g})$, where $V_0(X)$ is the
algebra of $C^\infty$ vector fields on $X$ having compact support
and $C_0^\infty(X;{\frak g})$ is the algebra of $C^\infty$
compactly supported  functions from $X$ into the Lie algebra
$\frak g$ of the  group $G$.  For each $(v,u)\in\frak a$, we
define the notion of a directional derivative of a function
$F\colon\OM\to\R$ along $(v,u)$, which is denoted by
$\nabla^\Omega_{(v,u)}F$. We obtain an explicit form of this
derivative on the special set $\FC$ of smooth cylinder functions
on $\OM$, which, in turn, motivates our definition of a tangent
bundle $T(\OM)$ of $\OM$, and of a gradient $\nabla^\Omega F$. We
note only that the tangent space $T_\om(\OM)$ to the marked
configuration space $\OM$ at a point $\om=(\gamma,s)\in\OM$ is
given by
\[ T_\om(\OM):= L^2(X\to T(X)\dotplus {\frak g};\gamma),\]
where $\dotplus$ means direct sum.

Next, we derive an integration by parts formula on $\OM$, that is,
we get an explicit formula for the dual operator $\div^\Om$ of the
gradient $\nabla^\Om$ on $\OM$. We prove that the probability
measures on $\OM$ for which $\nabla^\Om$ and $\div^\Om$ become
dual operators (with respect to $\la\cdot,\cdot\rangle_{T(\OM)}$)
are exactly the mixed marked Poisson measures
\[\mu_{\kappa,\tsigma}=\int_\Rp
\pi_{z\tsigma}\,\kappa(dz),\] where $\kappa$ is a probability
measure on $\Rp$ (with finite first moment) and $\pi_{z\tsigma}$
is the marked Poisson measure on $\OM$ with L\'evy measure
$z\tsigma$, $z\ge0$. This means that the mixed marked Poisson
measures are exactly the ``volume elements'' corresponding to our
differential geometry on $\OM$.

Thus, having identified the right volume elements on $\OM$, we
introduce for each measure $\mu_{\kappa,\tsigma}$ the first order
Sobolev space $H_0^{1,2}(\OM,\mu_{\kappa,\tsigma})$ by closing the
corresponding Dirichlet form
\[ {\cal E}^\Om_{\mu_{\kappa,\tsigma}}(F,G)=\int_{\OM}\la \nabla^\Om F,\nabla^\Om
G\ra_{T(\OM)}\, d\pi_{\kappa,\tsigma},\qquad F,G\in\FC,\] on
$L^2(\OM,\mu_{\kappa,\tsigma})$. Just as in the analysis on the
usual configuration space, this is the step where we really start
doing real infinite dimensional analysis. The corresponding
Dirichlet operator is denoted by $H_{\mu_{\kappa,\tsigma}}^\Om$;
it is a positive definite selfadjoint operator on
$L^2(\OM,\mu_{\kappa,\tsigma})$. The heat semigroup $\big( \exp(-t
H^\Om_{\mu_{\kappa,\tsigma}}) \big)_{t\ge0}$  generated by it is
calculated explicitly. The results on the ergodicity of this
semigroup are absolutely analogous to the corresponding results of
\cite{AKR}. Particularly, we have ergodicity if and only if
$\mu_{\kappa,\tsigma}=\pi_{z\tsigma}$ for some $z>0$, i.e.,
$\mu_{\kappa,\tsigma}$ is a (pure) marked Poisson measure.

We also clarify the relation between the intrinsic geometry on
$\OM$ we have constructed with another kind of extrinsic geometry
on $\OM$ which is based on fixing the marked Poisson measure
$\pi_\tsigma$ and considering the unitary isomorphism between
$L^2(\OM,\pi_\tsigma)$ and the corresponding Fock space
\[ {\cal F}(L^2(\XM;\tsigma))=\bigoplus_{n=0}^\infty\hat
L^2((\XM)^n,n!\,\tsigma^{\otimes n}),\] where $\hat
L^2((\XM)^n,n!\,\tsigma^{\otimes n})$ is the subspace of symmetric
functions from \linebreak $L^2((\XM)^n,n!\,\tsigma^{\otimes n})$.
Our main result here is to prove that $H^\Om_{\pi_\tsigma}$ is
unitarily equivalent (under the above isomorphism) to the second
quantization operator of the Dirichlet operator
$H_{\tsigma}^{\XM}$ on the $L^2(\XM;\tsigma)$
 space.

As a consequence of the results of this paper, we obtain a
representation on the marked Poisson space $L^2(\pi_\tsigma)$ not
only of the group $\frak A$, but also of its Lie algebra $\frak
a$. Let us remark that the groups of smooth (as well as measurable
and continuous) currents are classical objects in representation
theory, see e.g.\ \cite{Sergio,VGG1,GGV1,GGV2,VGG3,Ismagilov} and
references therein for different representations of these groups.
On the other hand, different representations of the group $\frak
A$ and its Lie algebra $\frak a$, in the special  case $G={\frak
g}=\R$, were constructed and studied by G.~Goldin et al.\
\cite{G1,G2,G3} from the point of view of nonrelativistic  quantum
mechanics.

Finally, we note that, in a way parallel to the work \cite{AKR2},
the results of the present paper can be generalized to the
interaction case where, instead of the Poisson measure
$\pi_{\tsigma}$, describing a system of free particles, one takes
its Gibbsian perturbation---more exactly, a marked Gibbs measure
on $\Omega_X^M$ of Ruelle type (see \cite{Kuna,KunaPhD}).

\section{Marked Poisson measures}

\subsection{Marked configuration space}

Let $X$ be a connected, oriented $C^\infty$ non-compact Riemannian
manifold. The configuration space $\Gamma_X$ over $X$ is defined
as the set of all locally finite subsets in $X$:
\[\Gamma_X:=\big\{\,\gamma\subset X\mid\#(\gamma\cap K)<\infty\ \text{for each
  compact } K\subset X\,\big\},\]
where $\#(\cdot)$ denotes the cardinality of a set. One can identify any
$\gamma\in\Gamma_X$ with the positive integer-valued Radon measure
\[\sum_{x\in\gamma}\eps_x\in{\cal M}(X),\]
where $\sum_{x\in\emptyset}\eps_x:=\text{zero measure}$ and ${\cal M}(X)$
 denotes the set of all positive  Radon measures on ${\cal B}(X)$.

Let also $M$ be a connected oriented $C^\infty$ (compact or non-compact)
Riemannian manifold.
The marked configuration space $\OM$ over $X$ with marks
from $M$ is defined as
\[ \OM:=\big\{\,\om=(\gamma,s)\mid \gamma\in\Gamma_X,\ s\in M^\gamma\,\big\},\]
where $M^\gamma$ stands for the set of all maps $\gamma\ni
x\mapsto m\in M$. Equivalently, we can define $\OM$ as the
collection of  subsets in $\XM$ having the following properties:
\[\OM=\left\{\om\subset X\times M\,\left|\, \begin{matrix}
 &\text{a)}\,\forall (x,m),(x',m')\in\om:(x,m)\ne(x',m')\Rightarrow x\ne x'\\
&\text{b)}\,\PR_X\om\in\Gamma_X\end{matrix}\right.\right\},\]
where $\PR_X$ denotes the projection of the Cartesian product of $X$
and $M$ onto $X$. Again, each $\om \in\OM$ can be identified with the measure
\[\sum_{(x,m)\in\om}\eps_{(x,m)}\in{\cal
  M}(X\times M).\]

It is worth noting that, for any bijection $\phi\colon \XM\to\XM$, the image of
the measure $\om(\cdot)$ under the mapping $\phi$, $(\phi^*\om)(\cdot)$,
coincides with $(\phi(\om))(\cdot)$, i.e.,
\[ (\phi^*\om)(\cdot)=(\phi(\om))(\cdot),\qquad \om\in\OM,\]
where $\phi(\om)=\{\phi(x,m)\mid(x,m)\in\om\}$ is the image of $\om$ as a
subset of $\XM$.

Let $\Bc(X)$ and $\Oc(X)$ denote the families of all Borel, resp.\ open subsets
of $X$ that have compact closure. Let also $\Bc(\XM)$ denote the family of all
Borel subsets of $\XM$ whose projection on $X$ belongs to $\Bc(X)$.

Denote by $C_{\text{0,b}}(\XM)$ the set of real-valued bounded continuous
functions $f$ on $\XM$ such that
$\supp f\in\Bc(\XM)$.
 As usually, we set for any $f\in C_{\text{0,b}}(\XM)$ and $\om\in\OM$
\[\la f,\om\ra=\int_\XM f(x,m)\,\om(dx,dm)=\sum_{(x,m)\in\om}f(x,m).\]
We note that, because of the definition of $\OM$, there are  only
a finite number of addends in the latter series.

Now, we are going to discuss the measurable structure of the space $\OM$. We
will use a ``localized'' description of the Borel $\sigma$-algebra ${\cal
  B}(\OM)$ over $\OM$.

For $\Lambda\in\Oc(X)$, define
\[\Om_\Lambda^{M}:=\big\{\,\om\in\OM\mid\PR_X\om\subset\Lambda\,\big\}\]
and for $n\in\Z_+=\{0,1,2,\dots\}$
\[ \Om_\Lambda^{M}(n):=\big\{\, \om\in\Om_\Lambda^{M}\mid\#(\om)=n\,\big\}.\]
It is obvious that
\[\Omega_\Lambda^{M}=\bigsqcup_{n=0}^\infty \Omega_\Lambda^M(n).\]

Let $\Lext:=\Lambda\times M$ (i.e., $\Lext$ is the set of all
``marked'' elements of  $\Lambda$) and let
\[{\widetilde \Lambda}^n_{\text{\ext}}:=\big\{
((x_1,m_1),\dots,(x_n,m_n))\in\Lambda_{\text{\ext}}^n\mid x_j\ne x_k\ \text{if }j\ne k\,\big\}.
\]
There is a bijection
\begin{equation}\label{one}
{\cal L}_\Lambda^{(n)}\colon {\widetilde\Lambda}_{\text{\ext}}^n/{\frak
  S}_n\mapsto
\Om_{\Lambda}^{M}(n)\end{equation}
given by
\[ {\cal L}_\Lambda^{(n)}\colon((x_1,m_1),\dots,(x_n,m_n))\mapsto
\{(x_1,m_1),\dots,(x_n,m_n)\}\in \Om_{\Lambda}^M(n),\] where
${\frak S}_n$ is the permutation group over $\{1,\dots,n\}$. On
$\Lext^n/{\frak S}_n$ one introduces the related metric
\begin{align*}
&
\delta\big[((x_1,m_1),\dots,(x_n,m_n)),((x_1',m_1'),\dots,(x_n',m_n'))\big]\\
&\qquad =\inf_{\sigma\in {\frak S}_n}
d^n\big[((x_1,m_1),\dots,(x_n,m_n)),((x'_{\sigma(1)},m'_{\sigma(1)}),\dots,(x'_{\sigma(n)},
m'_{\sigma(n)}))\big],
\end{align*}
where $d^n$ is the metric on $\Lext^n$ driven from the original metrics on $X$
and $M$. Then, ${\widetilde \Lambda}_{\text{\ext}}^n/{\frak S}_n$ becomes an open set in
$\Lext^n/{\frak S}_n$ and let ${\cal B}({\widetilde\Lambda}_{\text{\ext}}^n/{\frak S}_n)$ be the
trace $\sigma$-algebra on ${\widetilde \Lambda}_{\text{\ext}}^n/{\frak S}_n$ generated by ${\cal
  B}(\Lext^n/{\frak S}_n)$. Let then ${\cal B}(\Om_\Lambda^M(n))$ be the image
$\sigma$-algebra of ${\cal B}({\widetilde\Lambda}_{\text{\ext}}^n/{\frak S}_n)$ under the bijection
${\cal L}_\Lambda^{(n)}$ and let ${\cal B}(\Om_\Lambda^{M})$ be the
$\sigma$-algebra on $\Om_\Lambda^{M}$ generated by the usual topology of
(disjoint) union of topological spaces.

For any $\Lambda\in\Oc(X)$, there is a natural restriction map $
p_\Lambda\colon\OM\mapsto \Om_\Lambda^{M}$ defined by
\[\OM\ni\om\mapsto p_\Lambda (\om):=\om\cap\Lext\in\Om_\Lambda^M.\]
The topology on $\OM$ is defined as the weakest topology making all the mappings
$p_\Lambda$ continuous. The associated $\sigma$-algebra is denoted by ${\cal
  B}(\OM)$.

For each $B\in\Bc(\XM)$, we introduce a function
$N_B\colon\OM\to\Z_+=\{0,1,2,\dots\}$ such that
\begin{equation} N_B(\om):=\#(\om\cap B),\qquad \om\in\OM.\label{NB}\end{equation}
Then, it is not hard to see that ${\cal B}(\OM)$ is the smallest
$\sigma$-algebra on $\OM$ such that all the functions $N_B$ are
measurable.

\subsection{Marked Poisson measure}

In order to construct a marked Poisson measure, we fix:

(i) an intensity measure $\sigma$ on the underlying manifold $X$, which is
supposed to be  a   nonatomic Radon one,

(ii) a non-negative function
\[X\times{\cal B}(M)\ni (x,\Delta)\mapsto p(x,\Delta)\in\R_+
\]
such that, for $\sigma$-a.a.\ $x\in X$, $p(x,\cdot)$ is a finite
measure on $M$.

Now, we define a measure $\tsigma$ on $(\XM,{\cal B}(\XM))$ as follows:
\begin{equation}\label{1}
\tsigma(A)=\int_{A} p(x,dm)\,\sigma(dx),\qquad A\in{\cal B}(\XM) .
\end{equation}
We will suppose that the measure $\tsigma$ is infinite and  for any
$\Lambda\in\Bc(X)$
\begin{equation}\label{2}
\tsigma(\Lext)=\int_X {\bf
1}_\Lambda(x)p(x,M)\,\sigma(dx)<\infty,\end{equation} i.e.,
$p(x,M)\in L^1_{\mathrm loc}(\sigma)$.

Now, we wish to introduce a marked Poisson measure on $\OM$ (cf.\ e.g.\
\cite{kingman,Kall}).  To this end, we
take first the measure $\tsigma^{\otimes n}$ on $(X\times M)^n$, and for any
$\Lambda\in\Oc(X)$, $\tsigma^{\otimes n}$ can be considered as a finite measure
on $\Lambda_{\text{\ext}}^n$. Since $\sigma$ is nonatomic, we get
\[\tsigma^{\otimes n}
(\Lambda_{\text{\ext}}^n\setminus {\widetilde\Lambda}_{\text{\ext}}^n)=0
\]
and we can consider $\tsigma^{\otimes n}$ as a measure on
$(\widetilde\Lambda^n_{\text{\ext}}/ {\frak S}_n,{\cal
  B}(\widetilde\Lambda^n_{\text{\ext}}/{\frak S}_n))$
such that
$$\tsigma^{\otimes n}(\widetilde\Lambda^n_{\text{\ext}}/{\frak
S}_n)=\tsigma(\Lext)^n.$$

Denote by $\tsigma_{\Lambda,n}:=\tsigma^{\otimes n}\circ ({\cal
  L}_\Lambda^{(n)})^{-1}$ the image measure on $\Omega_\Lambda^M(n)$ under the
  bijection \eqref{one}. Then, we can define a measure
  $\lambda^\Lambda_{\tsigma}$ on $\Omega_\Lambda^{M}$ by
\[ \lambda_\tsigma^\Lambda:=\sum_{n=0}^\infty\frac1{n!}\,\tsigma_{\Lambda,n},\]
where $\tsigma_{\Lambda,0}:=\eps_{\emptyset}$ on
$\Omega_\Lambda^M(0)=\{\emptyset\}$. The measure $\lambda_\tsigma^\Lambda$ is
finite and
$\lambda^\Lambda_\tsigma(\Omega_\Lambda^M)=e^{\tsigma(\Lext)}$. Hence, the
measure
\[ \pi_\tsigma^\Lambda:=e^{-\tsigma(\Lext)}\lambda^\Lambda_{\tsigma}\]
is a probability measure on ${\cal B}(\Omega_\Lambda^M)$. It is not hard to
check the consistency property of the family
$\{\pi_\tsigma^\Lambda\mid\Lambda\in\Oc(X)\}$ and thus to obtain a unique
probability measure $\pi_\tsigma$ on ${\cal B}(\OM)$ such that
\[ \pi_\tsigma^\Lambda=p_\Lambda^*\pi_\tsigma,\qquad \Lambda\in\Oc(X).\]
This measure $\pi_\tsigma$ will be called a marked Poisson measure with L\'evy
measure $\tsigma$.

For any function $\fii\in C_{\text{0,b}}(\XM)$, it is easy to calculate the Laplace
transform of the measure $\pi_\tsigma$
\begin{equation}\ell_{\pi_\tsigma}(\fii):=\int_{\OM}e^{\la
    \fii,\om\ra}\,\pi_\tsigma(d\om)
=\exp\bigg(\int_{\XM}(e^{\fii(x,m)}-1)\,\tsigma(dx,dm)\bigg).\label{4}\end{equation}
\vspace{2mm}

\begin{exa}\rom{ Let $p(x,\cdot)\equiv \eps_m(\cdot)$, where $m$ is some
fixed point of $M$ and $x\in X$. Then,
$\tsigma=\sigma\otimes\eps_m$ and $\pi_\tsigma=\pi_\sigma$ is just
the Poisson measure on $(\Gamma_X,{\cal B}(\Gamma_X))$ with
intensity $\sigma$.}\end{exa}

\begin{exa}\rom{ Let $p(x,\cdot)\equiv \tau(\cdot)$, $x\in X$, where $\tau$
is a finite measure on $(M,{\cal B}(M))$. Now,
$\tsigma=\hat\sigma=\sigma\otimes\tau$ and $\pi_\tsigma$ coincides
with the marked Poisson measure under consideration in \cite{KSS}
(in the case where $M$ is a manifold). Notice that the choice of
$\tsigma=\hat\sigma$ as a product measure means a
position-independent marking, while the choice of a general
$\tsigma$ of the form \eqref{1} leads to a position-depending
marking.}\end{exa}

\section{Transformations of the marked Poisson measure}

\subsection{Group of transformations of the marked configuration space}

We are looking for a natural group $\Agot$ of transformations of $\OM$ such that
\begin{description}
\item[{\rm (i)}] $\pi_{\tsigma}$ is $\Agot$-quasiinvariant;

\item[{\rm(ii)}] $\Agot$ is big enough to reconstruct $\pi_{\tsigma}$ by the Radon--Nikodym
density $\dfrac{da^*\pi_\tsigma}{d\pi_\tsigma}$, where $a$ runs through $\Agot$.
\end{description}

Let us recall that in the work \cite{KSS} the group $\Diff$ was
taken as $\Agot$, just in the same way as in the case of the usual
Poisson measure \cite{AKR}. Here, $\Diff$ stands for the group of
diffeomorphisms of $X$ with compact support, i.e., each
$\psi\in\Diff$ is a diffeomorphism of $X$ that is equal to the
identity outside a compact set (depending on $\psi$). The group
$\Diff$ satisfies (i). However, unlike the case of the Poisson
measure, the condition (ii) is not  satisfied, because, for
example, in the case where $\tsigma=\sigma\otimes\tau$, there is
no information about the measure $\tau$ that is contained in
$\dfrac{d\psi^*\pi_\tsigma}{d\pi_\tsigma}$, see \cite{KSS}.
Therefore, just as in the case of \cite{KLU1}, we need a proper
extension of the group $\Diff$.

In what follows, we will suppose that  $M$ is a homogeneous space of a Lie group
$G$ (see e.g.\ \cite{Bo}). Let us recall that this means the existence of a
$C^\infty$ mapping $\theta\colon G\times M\to M$ satisfying the following
conditions:
\begin{description}
\item[{\rm (i)}] If $e$ is the unity element of the group $G$, then
\[ \theta(e,m)=m\qquad \text{for all }m\in M;\]
\item[{\rm (ii)}]
 If $g_1,g_2\in G$, then
\[ \theta (g_1,\theta(g_2,m))=\theta(g_1g_2,m)\qquad \text{for all }m\in M;\]
\item[{\rm (iii)}] For  arbitrary $m_1,m_2\in M$, there exists $g\in G$ such that
$\theta(g,m_1)=m_2$.
\end{description}

For any $g\in G$, we will denote by $\theta_g\colon M\to M$ the mapping
given by $\theta_g(m)\colon =\theta(g,m)$; then $\theta_g$ defines a
diffeomorphism of $M$.

Let us fix an arbitrary point $m_0\in M$ and let $H$ be the isotropy group of $M$:
\[ H:=\big\{\,g\in G\mid \theta_g(m_0)=m_0\,\big\}.\]
Then, the homogeneous space $M$ can always be identified with the
factor space $G/H$ (endowed with the unique corresponding
$C^\infty$ manifold structure), i.e., $M=G/H$.

Let us consider the group of {\it smooth  currents}, i.e., all
$C^\infty$ mappings $X\ni x\mapsto\eta(x)\in G$, which are equal
to $e$ outside a compact set (depending on $\eta$). A
multiplication $\eta_1\eta_2$ in this group is defined as the
pointwise multiplication of the mappings $\eta_1$ and $\eta_2$. In
the representation theory this group is denoted by $G^X$, or
$C^\infty_0(X;G)$.

The group $\Diff$ acts in $G^X$  by automorphisms: for each $\psi\in \Diff$,
\[G^X\ni\eta\overset{\alpha}{\mapsto}\alpha(\psi)\eta:=\eta\circ\psi^{-1}\in G^X.\]
 Thus, we can endow the Cartesian product of $\Diff$ and $G^X$ with the
 following multiplication: for $a_1=(\psi_1,\eta_1)$, $a_2=(\psi_2,\eta_2)$ from
 $\Diff\times G^X$
\[ a_1a_2=(\psi_1\circ\psi_2,\eta_1(\eta_2\circ\psi_1^{-1}))\]
and obtain a semidirect product
\[ \Diff\underset{\alpha}{\times}G^X=:\Agot\]
of the groups $\Diff$ and $G^X$.

The group $\Agot$ acts in $\XM$ in the following way: for any
$a=(\psi,\eta)\in\Agot$
\begin{equation}\label{3}
\XM\ni(x,m)\mapsto
a(x,m)=(\psi(x),\theta(\eta(\psi(x)),m))\in\XM.\end{equation} If
id denotes the identity diffeomorphism of $X$ and ${\pmb{e}}$ is
the function identically equal to $e$ on $X$, then we will just
identify $\psi$ with $(\psi,{\pmb{e}})$ and $\eta$ with
$(\operatorname{id},\eta)$. The action \eqref{3} of an arbitrary
$a=(\psi,\eta)$ can be represented as
\[(x,m)\mapsto a(x,m)=\eta\psi(x,m),\]
where
\begin{align*}
\psi(x,m)&=(\psi(x),m),\\
\eta(x,m)&=(x,\theta(\eta(x),m)).\end{align*}

For any $a=(\psi,\eta)\in\Agot$, denote $K_a:=K_\psi\cup K_\eta$, where
$K_\psi$ and $K_\eta$ are the minimal closed sets in $X$ outside of which
$\psi=\operatorname{id}$ and $\eta=\pmb {e}$, respectively. Evidently,
$K_a\in{\cal B}_{\mathrm c}(X)$,
\[a(K_a)_{\mathrm \ext}=(K_a)_{\mathrm \ext},\]
and $a$ is the identity transformation outside $(K_a)_{\mathrm \ext}$.

Now, let us recall some known facts concerning quasiinvariant measures on
homogeneous spaces (see e.g.\ \cite{Wa,VilKlimyk}).

\begin{th}\label{op} Suppose $G$ is a Lie group and $H$ its subgroup\rom, and
  let $dg$\rom, $\delta_G$ and $dh$\rom, $\delta_H$ be fixed  Haar measures and
  modular functions on $G$ and $H$\rom, respectively\rom. Then\rom:
\begin{description}

\item{\rom{(i)}} for every measure $\mu$ on $G/H$ that is quasiinvariant with
  respect to the action of $G$ on $G/H$\rom, there exists a measurable positive
  function $\xi$ on $G$ verifying
\begin{equation}\label{una}
\xi(gh)=\frac{\delta_H(h)}{\delta_G(h)}\,\xi(g),\qquad g\in G,\, h\in
  H,\end{equation}
and
\begin{equation}
\label{duo}
\int _G f(g)\xi(g)\, dg=\int_{G/H}\mu(d\, gH)\int_H f(gh)\,dh,\qquad f\in
C_0(G),\end{equation}
where $C_0(G)$ denotes the set of continuous functions on $G$ with compact
support\rom; for each $g\in G$ the Radon--Nikodym density is given by
\[ p_g^\mu(\tilde gH):=\frac{dg^*\mu}{d\mu}(\tilde gH)=\frac{\xi(g^{-1}\tilde
  g)}
{\xi(\tilde g)},\qquad \tilde gH\in G/H;\]

\item{\rom{(ii)}} there exists a quasiinvariant measure $\lambda$ on $G/H$ such
  that
the function $$p^\lambda(g,\tilde gH):=p_g^\lambda(\tilde gH)$$ is differentiable
  on $G\times G/H$\rom.
\end{description}
\end{th}

\begin{rem}\rom{
We recall that the modular function $\delta_G(\cdot)$ of a Lie
group $G$ is defined from the equality $r^*_{\tilde
g}\,dg=\delta_G(\tilde g)\,dg$, where $dg$ is the Haar measure on
$G$ (i.e., a fixed left-invariant measure on $G$) and $r_g$
denotes the right translation on $G$, i.e., $\tilde g\mapsto
r_g\tilde g=g\tilde g$. }\end{rem}

We fix the measure $\lambda$ on $M=G/H$ from Theorem~\ref{op}, (ii). As easily
seen from Theorem~\ref{op} (i),  any quasiinvariant measure on $M$ in
equivalent to $\lambda$.

\begin{rem}\rom{
If $H=\{e\}$, i.e., $M=G$, then we can choose $\lambda$ to be the
Haar measure $dg$ on $G$. Moreover, if $\delta_G(h)=\delta_H(h)$
for all $h\in H$ (and only in this case) there exists a  $\lambda$
being invariant with respect to the action of $G$ on $M$. The
latter condition holds automatically if  $G$ is unimodular, that
is, $\delta_G(g)\equiv 1$ for all $g\in G$. This, in turn, holds
for all compact and simple Lie groups. }\end{rem}

In what follows, we will suppose that the measure $\sigma$ is equivalent to the
Riemannian volume $\nu$ on $X$: $\sigma(dx)=\rho(x)\,\nu(dx)$ with $\rho>0$ $\nu$-a.s.,
and that for $\nu$-a.a.\ $x\in X$ $p(x,\cdot)$ is equivalent to the measure
$\lambda$:
\[ p(x,dm)=p(x,m)\,\lambda(dm)\quad\text{with }p(x,m)>0\ \lambda\text{-a.a.}\ m\in
M.\]
Thus, the measure $\tsigma$ can be written in the form
\[\tsigma(dx,dm)=\rho(x)p(x,m)\,\nu(dx)\,\lambda(dm).\]
The condition $\tsigma(\Lext)<\infty$, $\Lambda\in\Bc(X)$, implies that the
function
\[ q(x,m):=\rho(x)p(x,m)\]
satisfies
\begin{equation}\label{c1}
q^{1/2}\in L^2_{\text{loc}}(X;\nu)\otimes L^2(M;\lambda).\end{equation}

Noting that
$$a^{-1}(x,m)=(\psi,\eta)^{-1}(x,m)=(\psi^{-1}(x),\theta(\eta^{-1}(x),m)),$$
we easily deduce the following

\begin{prop}\label{prop1} The measure
$\tsigma$ is $\Agot$-quasiinvariant and for any $a=(\psi,\eta)\in\Agot$
the Radon--Nikodym density is given by
\[
\left\{
\begin{aligned}
\text{}& p_a^{\tsigma}(x,m):=\frac{d(a^*\tsigma)}{d\tsigma}(x,m)=
\frac{q(\psi^{-1}(x),\theta(\eta^{-1}(x),m)
)}{q(x,m)}p^\lambda(\eta(x),m)\,J_\nu^\psi(x),\\
 \text{}&\text{\rom{if}
  }(x,m)\in\{0<q(x,m)<\infty\}\cap\{0<q(\psi^{-1}(x),\theta(\eta^{-1}(x),m)
)<\infty\},\\
\text{}&p_a^\tsigma(x,m)=1, \qquad \text{\rom{otherwise}},
\end{aligned}
\right.\]
where $J_\nu^\psi$ is the Jacobian determinant of $\psi$ \rom(w\rom.r\rom.t\rom.\
the Riemannian volume $\nu$\rom{).}
\end{prop}

We give two examples of the above construction, which are
important from the point of view of the marked configuration space
analysis. We refer the reader to e.g.\ \cite{VilKlimyk, Wa} for
further examples.

\begin{exa}\rom{Let $G=\R_+$ be the dilation group (e.g.\ \cite{Goldin}), i.e., the
multiplication in this group is given by the usual multiplication
of numbers. As a homogeneous space $M$ we take  $G$ itself, by
identifying the action of the group with the multiplication in it.
As a quasiinvariant measure $\lambda$ on $M$ we can take the
restriction to $\R_+$ of the Lebesgue measure on $\R$.

The analysis and geometry on the marked configuration space
$\Omega_X^{\R_+}$ were studied in our previous work \cite{KLU1}.
Here we only mention that the choice $M=\R_+$ leads (via a natural
isomorphism) to the class of compound Poisson measures. In other
words, each mark $s_x\in\R_+$ corresponding to $x\in X$
 describes the charge of the measure $$\omega=(\gamma,s)=\sum_{x\in
 X}s_x\eps_x\in{\cal M}(X)$$ at the point $x$ (or, in the case where
 $X=\R$,
 the value of the jump of the process at $x$).
}\end{exa}

\begin{exa}\rom{Let $G=O(d+1)$ be the $(d+1)$-dimensional
orthogonal group and let $M=S^d$ be the $d$-dimensional unit
sphere in $\R^{d+1}$ with the natural action of the group $O(d+1)$
on $S^d$, see e.g.\ \cite{GMS, VilKlimyk, Wa}. As $\lambda$ we
take the surface measure on $S^d$, which is invariant w.r.t.\ the
action of $O(d+1)$.  From the point of view of statistical
mechanics, a mark $s_x\in S^d$  describes in this example the spin
of the particle at the point $x$. }\end{exa}

\subsection{$\Agot$-quasiinvariance of the marked Poisson measure}

Any $a\in\Agot$ defines by \eqref{3} a transformation of $X\times M$, and,
consequently, $a$ has the following ``lifting'' from $\XM$ to $\OM$:
\begin{equation}\label{kk} \OM\ni\om\mapsto a(\om)=\big\{\,
a(x,m)\mid (x,m)\in\om\,\big\}\in\OM.\end{equation}
(Note that, for a given $\om\in\OM$, $a(\om)$ indeed belongs to $\OM$ and
coincides with $\om$ for all but a finite number of points.)
The mapping \eqref{kk} is obviously measurable and we can define the image
$a^*\pi_\tsigma$ as usually.
The following proposition is an analog of a corresponding fact about Poisson
measures.

\begin{prop}\label{prop2} For any $a\in\Agot$\rom, we have
\[ a^*\pi_\tsigma=\pi_{a^*\tsigma}.\]
\end{prop}

\noindent{\it Proof}. The proof is the same as for the usual  Poisson
measure $\pi_\sigma$ with intensity $\sigma$ and $\psi\in\Diff$ (e.g.,
\cite{AKR}), one has just to calculate the Laplace transform of the measure
$a^*\pi_\tsigma$ for any $f\in C_{\text{0,b}}(\XM)$ and to use the formula
\eqref{4}.\quad$\blacksquare$

\begin{prop}\label{prop3} The marked Poisson measure $\pi_{\tsigma}$ is
  quasiinvariant w.r.t.\ the group $\Agot$, and for any $a\in\Agot$ we have
\begin{equation}\label{rr}
\frac{d(a^*\pi_\tsigma)}{d\pi_\tsigma}(\om)=\prod_{(x,m)\in\om}p_a^\tsigma
(x,m)
.\end{equation}
\end{prop}

\noindent{\it Proof}. The result follows from  Skorokhod theorem
on absolute continuity of Poisson measures (see, e.g.,
\cite{Sk,Ta}). \quad$\blacksquare$

\begin{rem}\rom{
Notice that only a finite (depending on $\omega$) number of factors in the product on the right hand side of \eqref{rr}
are not equal to one.
}\end{rem}

\section{The differential geometry of marked configuration spaces}

\subsection{The tangent bundle of $\Omega^{M}_{X}$}

\newcommand{\agot}{{\frak a}}

Let us denote by $V_0(X)$ the set of $C^\infty$ vector fields on
$X$ (i.e., smooth sections of $T(X)$) that have compact support.
Let $\got$ denote the Lie algebra of $G$ and let
$C_0^\infty(X;\got)$ stand for the set of all $C^\infty$ mappings
of $X$ into $\got$ that have compact support. Then $$\agot:
=V_0(X)\times C_0^\infty(X;\got)$$ can be thought of as a Lie
algebra corresponding to the Lie group $\Agot$.
 More precisely, for any fixed
$v\in V_0(X)$ and for any $x\in X$, the curve
\[\R\ni t\mapsto\psi_t^v(x)\in X\]
is defined as the solution of the following Cauchy problem
\begin{equation}\label{beta}\left\{
\begin{aligned}
\frac{d}{dt}\psi_t^v(x)&=v(\psi_t^v(x)),\\
\psi_0^v(x)&=x.\end{aligned}\right.
\end{equation}
Then, the mappings $\{\psi_t^v,\, t\in\R\}$ form a one-parameter subgroup of
diffeomorphisms in $\Diff$ (see, e.g., \cite{Bo}):
\begin{align*}
1)&\forall t\in\R\quad \psi_t^v\in\Diff,\\
2)&\forall t_1,t_2\in\R\quad \psi_{t_1}^v\circ\psi_{t_2}^v=\psi^v_{t_1+t_2}.
\end{align*}
Next, for each function $u\in C^\infty_0(X;\got)$, $x\in X$, and
$t\in\R$, we set $ \eta_t^u(x):=\exp(t u(x))$, where $\got\ni
Y\mapsto \exp Y\in G$ is the exponential mapping (see, e.g.,
\cite{Wa}). Hence, for a fixed $x\in X$, $\{\eta_t^u(x),\,
t\in\R\}$ is  a one-parameter subgroup of $G$ and
\begin{equation}\label{3.1prime}
\begin{gathered}
\eta_0^u(x)=e,\\
\frac{d}{dt}\,\eta_t^u(x)\big|_{t=0}=u(x).\end{gathered}\end{equation}
Let us recall a fundamental theorem in the theory of Lie groups.

\begin{th} There exists a neighborhood $U$ of the zero in $\got$ and a neighborhood
$O$ of the unit element $e$ in $G$ such that $\exp\colon U\to O$ is an analytic
diffeomorphism.\end{th}

From this theorem, we conclude that, for each fixed $u\in C_0^\infty(X;\got)
$, there exists $\eps>0$ such that for any $t\in(-\eps,\eps)$ the mapping
$X\ni x\mapsto\eta_t^u(x)\in G$ belongs to $G^X$, which yields, in turn,
that $\eta_t^u\in G^X$ for all $t\in\R$, and moreover $\eta_t^u$ is a
one-parameter subgroup of $G^X$.

Thus, for an arbitrary $(v,u)\in \agot$, we can consider the curve
$\{(\psi_t^v,\eta_t^u),\, t\in\R\}$ in $\Agot$. Hence, to any
$\om\in\OM$ there corresponds the following curve in $\OM$:
$$\R\ni t\mapsto (\psi_t^v,\eta_t^u)\om\in\OM.$$ Define now for a
function $F\colon\OM\to\R$ the directional derivative of $F$ along
$(v,u)$ as
\[ (\nabla^\Om_{(v,u)}F)(\om):=\frac
d{dt}F((\psi_t^v,\eta_t^u)\om)\big|_{t=0},\]
provided the right hand side exists.
We will also denote by $\nabla^\Om_v$ and
$\nabla^\Om_u$ the directional derivatives along $(v,0)$ and $(0,u)$,
respectively.

Absolutely analogously, one defines for a function $\fii\colon \XM\to\R$ the
directional derivative of $\fii$ along $(v,u)$:
\begin{equation}\label{small}
(\nabla^\XM_{(v,u)}\fii)(x,m)=\frac{d}{dt}\fii((\psi_t^v,\eta_t^u)(x,m))
\big|_{t=0}.\end{equation}
Then, for a continuously differentiable function $\fii$, we have from
\eqref{3}, \eqref{beta}, \eqref{3.1prime}, and \eqref{small}
\begin{gather}
(\nabla^\XM_{(v,u)}\fii)(x,m)=
\frac{d}{dt}\,\fii((\psi_t^v(x),\theta
(\eta_t^u(\psi_t^v(x)),m)\big|_{t=0}\notag\\
=\frac d{dt}\fii(\psi_t^v(x),m)\big|_{t=0}+\frac d{dt}\fii(x,\theta(\eta_t^u(x)
,m))\big|_{t=0}\notag\\
\mbox{}+ \frac{d}{dt}\fii(x,\theta(\eta_0^u(\psi_t^v(x)),m))\big|_{t=0}\notag\\
=\la\nabla^X\fii (x,m),v(x)\ra_{T_x(X)}+\la\nabla^G\fii(x,\theta(e,m)),u(x)\ra
_{\got}\notag\\
=\la\nabla^\XM\fii(x,m),(v(x),u(x))\ra_{T_{(x,m)}(\XM)}.
\label{9}
\end{gather}
Here, $ T_{(x,m)}(\XM):= T_x(X)\dotplus\got$ and
$\nabla^\XM:=(\nabla^X,\nablat^{M})$, where $\nabla^X$
 denotes the gradient on $X$ and
\begin{equation}\label{bum}\begin{gathered}
\nablat^{M}f(m)=\nabla^G \hat f(e,m),\\ \hat
f(g,m):=f(\theta(g,m)),\qquad g\in G,\, m\in M,
\end{gathered}\end{equation}
$\nabla^G$ being the gradient on $G$.

\begin{rem}\label{brandnewremark}\rom{
Notice that upon \eqref{bum} we have, for a fixed $u\in\got$,
\begin{align}
\la\nablat^M f(m),u\ra_\got&=\la\nabla^G
f(\theta(e,m)),u\ra_\got\notag\\ &=\frac
d{dt}f(\theta(e^{tu},m))\big|_{t=0}\notag\\ &=\la\nabla^M f(m),
(Ru)(m)\ra_{T_m(M)},\label{ghasdf}\end{align} where $\nabla^M$
denotes the usual gradient on $M$, and the vector field $Ru$ on
$M$ is given by
\begin{equation}\label{oqhdqpo}M\ni m\mapsto (Ru)(m):=\frac
d{dt}\theta(e^{tu},m)\big|_{t=0}.\end{equation}  }\end{rem}

Let us introduce a special class of ``nice functions'' on $\OM$.
 Denote by $\frak D$ the set of all  $C^\infty$-functions $\fii$ on
 $\XM$  such that  the support of $\fii$ is in  $\Bc(\XM)$, and
 $\fii$ and all its
$\nabla^\XM$
derivatives are
 bounded.
 Next, let
$C_{\text{b}}^\infty(\R^N)$ stand for the space of all
$C^\infty$-functions on $\R^N$ which together with all their derivatives are
bounded. Then, we can introduce $\FC$
as the set of all functions $F\colon\OM\mapsto\R$ of the form
\begin{equation}\label{5}
F(\om)=g_F(\la\vph_1,\om\ra,\dots,\la\vph_N,\om\ra),\qquad
\om\in\OM,\end{equation}
where $\vph_1,\dots,\fii_N\in{\frak D}$ and $g_F\in C_{\mathrm b}^\infty(\R^N)$
(compare with \cite{AKR}). $\FC$ will be called the set of smooth cylinder
functions on $\OM$.

For any  $F\in\FC$ of the form \eqref{5} and a  given $(v,u)\in \agot$, we have,
just as in \cite{AKR},
\begin{align*}F((\psi_t^v,\eta_t^u)\om)&=g_F(\la\fii_1,(\psi_t^v,\eta_t^u)
\om\ra,\dots,
\la\fii_N,(\psi_t^v,\eta_t^u)\om\ra)\\
&=g_F(\la\fii_1\circ(\psi_t^v,\eta_t^u),\om\ra,\dots,\la\fii_N\circ
(\psi_t^v,\eta_t^u),\om\ra),\end{align*}
and therefore
\begin{equation}\label{8}
(\nabla_{(v,u)}^\Omega F)(\om)=\sum_{j=1}^N\frac
{\di g_F}{\di r_j}(\la\fii_1,\om\ra,\dots,\la\fii_N,\om\ra)
\la\nabla_{(v,u)}^\XM\fii_j,\om\ra.\end{equation}
In particular, we conclude from \eqref{8} that
\begin{equation}\label{10}
\nabla^\Om_{(v,u)}=\nabla^\Om_v+\nabla_u^\Om.
\end{equation}

The expression of $\nabla_{(v,a)}^{\Omega}$ on smooth cylinder functions
motivates the following definition.

\begin{define}\label{def0}\rom{ The tangent space $T_\om\big(\OM\big)$
to the marked
    configuration space $\OM$ at a point $\om=(\gamma,s)\in\OM$ is defined
    as the Hilbert space
\begin{align*}
T_\om\big(\OM\big):&=L^2(X\to T(X)\dotplus\got;\gamma)\\
&=L^2(X\to T(X);\gamma)\oplus
L^2(X\to\got;\gamma)\\&=\bigoplus_{x\in\gamma}\big[T_x(X)\oplus
\got\big]
\end{align*}
with  scalar product \begin{align}\la
V_\om^1,V_\om^2\ra_{T_\om\left(\OM\right)} &=\int_X\big(\la
V_\om^1(x)_{T_x(X)},V_\om^2(x)_{T_x(X)}\ra_{T_x(X)}+ \la
V_\om^1(x)_{\got},
V_\om^2(x)_{\got}\ra_{\got}\big)\,\gamma(dx)\notag\\
&=\sum_{x\in\gamma}\big(\la
V_\om^1(x)_{T_x(X)},V_\om^2(x)_{T_x(X)}\ra_{T_x(X)}+ \la
V_\om^1(x)_{\got}, V_\om^2(x)_{\got}\ra_{\got}\big),\label{9.1}
\end{align} where $V_\om^1,V_\om^2\in T_\om\big(\OM\big)$ and
$V_\om(x)_{T_x(X)}$ and $V_\om(x)_{\got}$ denote the projection of
\linebreak $V_\om(x)\in T_x(X)\dotplus\got$ onto $T_x(X)$ and
$\got$, respectively. (Notice that the tangent space
$T_\om\big(\OM\big)$ depends only on the $\gamma$ coordinate of
$\omega$.) The corresponding tangent bundle is
\[ T\big(\OM\big)=\bigcup_{\om\in\OM}T_\om\big(\OM\big).\]
}\end{define}

As usually in Riemannian geometry, having directional derivatives
and a Hilbert space as a tangent space, we can introduce a
gradient.

\begin{define}\label{def1}
\rom{ We define the intrinsic gradient $\nabla^ \Omega$ of a
function $F\colon \OM\to\R$ as the mapping
\[\OM\ni\om\mapsto(\nabla^\Omega F)(\om)\in T_\om\big(\OM\big)\]
such that, for any $(v,u)\in\agot$,
\[(\nabla^\Omega_{(v,u)}F)(\om)=\la(\nabla^\Omega F)(\om),(v,u)
\ra_{T_\om\left(\OM\right)}.\]
}\end{define}

By \eqref{8} and \eqref{9} we have, for an arbitrary $F\in\FC$ of the form
\eqref{5} and each $\om=(\gamma,s)\in\OM$,
\begin{equation}\label{gradient}(\nabla^\Omega F)(\om;x)=\sum_{j=1}^N
\frac{\di g_F}{\di
  r_j}(\la\fii_1,\om\ra,\dots,\la\fii_N,\om\ra)
\nabla^\XM\fii_j(x,s_x),\qquad x\in\gamma.\end{equation}

\subsection{Integration by parts and divergence on the marked Poisson space}

Let the marked configuration space $\OM$ be equipped with the
marked Poisson measure $\pi_{\tsigma}$. We strengthen the
condition \eqref{c1} by demanding that
\begin{equation}\label{condition} q^{1/2}\in H_0^{1,2}(\XM).\end{equation}
Here, $H_0^{1,2}(\XM)$ denotes the local Sobolev space of order 1
constructed with respect
to the gradient $\nabla^\XM$ in the space $L^2_{\mathrm loc}(X;\nu)\otimes
L^2(M;\lambda)$, i.e., $H_0^{1,2}(\XM)$ consists of functions $f$ defined on
$\XM$ such that, for any set $A\in{\cal B}_{\mathrm c}(\XM)$, the
restriction of
$f$ to $A$ coincides with the restriction to $A$ of some function $\fii$ from the
Sobolev space $H^{1,2}(\XM)$ constructed as the closure of $\frak D$ with
respect to the norm
\[ \|\fii\|^2_{1,2}:=\int_{\XM}\Big( |\nabla^X\fii(x,m)|^2_{T_x(X)}+
|\nablat^M\fii(x,m)|_{\got}^2
+|\fii(x,s)|^2\Big)\nu(dx)\,\lambda(dm).\]

Additionally, we will suppose that, for each $\Lambda\in\Bc(X)$,
\begin{equation}\label{neues}
|\nabla^G p^\lambda(e,\cdot)|_{\got}\in L^1(\Lext,\tsigma),\end{equation}
where, as before,
$$p^\lambda(g,m)=\frac{dg^*\lambda}{d\lambda}(m).$$

The set $\FC$ is a dense subset in the space
\[ L^2(\OM,{\cal B}(\OM),\pi_{\tsigma})=:L^2(\pi_{\tsigma}).\]
For any $(v,u)\in\agot$, we have a differential operator
 in $L^2(\pi_\tsigma)$ on the domain \linebreak$\FC$ given by
\[\FC\ni F\mapsto\nabla_{(v,u)}^\Omega F\in L^2(\pi_\tsigma).\]
Our aim now is to compute the adjoint operator
$\nabla_{(v,u)}^{\Omega\,*}$ in $L^2(\pi_\tsigma)$. This
corresponds, of course, to the deriving of an integration by parts
formula with respect to the measure $\pi_\tsigma$.

But first we present the corresponding formula on $\XM$.

\begin{define}\rom{
For any $(v,u)\in\agot$, the
logarithmic derivative of the measure $\tsigma$
along $(v,u)$ is defined as the following function on $\XM$:
\[\beta_{(v,u)}^\tsigma:=\beta_v^\tsigma+\beta_u^\tsigma\]
with
$$
\beta_v^\tsigma(x,m)=\bigg\la\frac{\nabla^Xq(x,m)}{q(x,m)},
v(x)\bigg\ra_
{T_x(X)}+
\div^Xv(x),
$$
$\div^X=\div^X_\nu$ being the divergence on $X$ w\rom.r\rom.t\rom.\
$\nu$\rom, and
\[\beta_u^\tsigma(x,m)
=\bigg\la\frac{\nablat^M q(x,m)}{q(x,m)},u(x)\bigg\ra_{\got}+\la\nabla^G p^\lambda
(e,m),-u(x)\ra_{\got}
.\]
}\end{define}

Upon \eqref{condition}, we conclude that, for each $(v,u)\in\agot$,
the function
$\nabla_{(v,u)}^\XM\log q$ is quadratically integrable with respect to the
measure $\tsigma$, and  therefore, since the support of $\nabla_{(v,u)}^\XM
\log q$ belongs to $\Bc(\XM)$, this function is from $L^1(\XM,\tsigma)$.
Thus, in virtue of the condition \eqref{neues}, we get the inclusion $
\beta_{(v,u)}^\tsigma\in L^1(\XM,\tsigma)$.

By using standard arguments, one shows the following

\begin{lem}[Integration by parts formula on $\XM$]
For all $\fii_1$\rom, $\fii_2\in{\frak D}$\rom, we have
\begin{align*}
&\int_\XM(\nabla^\XM_{(v,u)}\fii_1)(x,m)\fii_2(x,m)\,\tsigma(dx,dm)=\\
&\qquad=-\int_\XM\fii_1(x,m)(\nabla^\XM
_{(v,u)}\fii_2)(x,m)\,\tsigma(dx,dm)
\\&\qquad\quad
\text{}-\int_\XM\fii_1(x,s)\fii_2(x,s)\beta^\tsigma_{(v,u)}(x,m)\,\tsigma(dx,dm).\end{align*}
 \end{lem}

\begin{rem}
\label{kiii} \rom{ The function $\la\nabla^G
p^\lambda(e,m),-u(x)\ra_{\got}$, which appears in the definition
of $\beta_u^{\tsigma}$ is, for each fixed $x\in X$, the divergence
on $M$ with respect to the measure $\lambda$ of the vector field
$Ru(x)$ on $M$ defined by \eqref{oqhdqpo}, see
Remark~\ref{brandnewremark}. Indeed, for any  $u\in\got$ and for
an arbitrary $f$ from $C_0^\infty(M)$---the space of all
$C^\infty$ functions on $M$ with compact   support, we have
\begin{align*}
\int_M\nablat_u^M f(m)\,\lambda(dm)&= \int_M\la \nabla^M
f(m),(Ru)(m)\ra_{T_m(M)} \,\lambda(dm)\\ &=\int_M\frac
d{dt}f(\theta(\exp(tu),m))\big|_{t=0}\,\lambda(dm)\\
&=\int_Mf(m)\,\frac
d{dt}p^\lambda(\exp(tu),m)\big|_{t=0}\,\lambda(dm)\\ &=\int_M
f(m)\la\nabla^G p^\lambda(e,m),u\ra_{\got}\,\lambda(dm).
\end{align*}
}\end{rem}

\begin{define}\rom{
For any $(v,u)\in\agot$, the logarithmic derivative of the marked
Poisson measure
$\pi_\tsigma$ along $(v,u)$ is defined as the following function on $\OM$:
\begin{equation}\label{16}\OM\ni\om\mapsto
B^{\pi_\tsigma}_{(v,u)}(\om):=\la\beta_{(v,u)}^\tsigma,\om\ra.\end{equation}
}\end{define}

A motivation for this definition is given by the following theorem.

\begin{th}
[Integration by parts formula]
For all $F_1,F_2\in\FC$ and each $(v,u)\in\agot$\rom, we have
\begin{align}
\int_{\OM}(\nabla_{(v,u)}^\Omega F_1)(\om)F_2(\om)\,\pi_\tsigma(d\om)&=-\int_{\OM}
F_1(\om)(\nabla_{(v,u)}^\Omega F_2)(\om)\,\pi_\tsigma(d\om)\notag\\
&\quad-\int_{\OM} F_1(\om) F_2(\om) B_{(v,u)}^{\pi_\tsigma}(\om)\,
\pi_\tsigma(d\om),
\label{11}\end{align}
or
\begin{equation}
\label{12}
\nabla_{(v,u)}^{\Om\,*}=-\nabla_{(v,u)}^\Omega-B_{(v,u)}^{\pi_{\tsigma}}(\om)\end{equation}
as an operator equality on the domain $\FC$ in $L^2(\pi_\tsigma)$\rom.
\label{thIbP}\end{th}

\noindent{\it Proof}. Because of \eqref{10}, the formula \eqref{12} will be
proved if we prove it first for the operator $\nabla^\Om_v$, i.e., when $u(x)
\equiv
0$, and then for the operator $\nabla^\Om_u$, i.e., when $v(x)=0\in T_x(X)$
for all
$x\in X$. We present below only the proof for $\nabla^\Om_u$, since the
proof for
$\nabla^\Om_v$ is basically the same as that of the integration by parts
formula in case of Poisson measures \cite{AKR}.


By Proposition~\ref{prop2}, we have for all $u\in
C^\infty_0(X;\got)$
\[\int_{\OM} F_1(\eta_t^u(\om))F_2(\om)\,\pi_\tsigma(d\om)=\int_{\OM}
F_1(\om)F_2(\eta_{-t}^u(\om))\,
\pi_{\eta_t^{u\,*}\tsigma}(d\om).\]
Differentiating this equation with respect to $t$, interchanging $d/dt$ with the
integrals and setting $t=0$, the l.h.s.\ becomes the l.h.s.\ of \eqref{11}. To see
that the r.h.s.\ then also coincides with the r.h.s.\ of \eqref{11}, we note
that
\[\frac d{dt}F_2(\eta_{-t}^u(\om))\big|_{t=0}=-(\nabla_u^\Omega F_2)(\om),\]
and by Proposition~\ref{prop3}
\begin{gather*}
\frac
d{dt}\bigg[\frac{d\pi_{\eta_t^{u\,*}\tsigma}}{d\pi_\tsigma}(\om)\bigg]
\big|_{t=0}=
\sum_{(x,m)\in\om}\frac d{dt}\, p_{\eta_t^u}^\tsigma (x,m)
\bigg|_{t=0}\\
=-\la\beta_u^\tsigma,\om\ra
=-B_u^{\pi_\tsigma}(\om).\quad\blacksquare
\end{gather*}

\begin{define}\label{def2}\rom{
For a vector field
\[V\colon\OM\ni\om\mapsto V_\om\in T_\om(\OM),\]
the divergence $\div_{\pi_\tsigma}^\Omega V$ is defined via the duality relation
\[ \int_{\OM}\la V_\om,\nabla^\Omega
F(\om)\ra_{T_\om\big(\OM\big)}\pi_\tsigma(d\om)=-\int_{\OM}F(\om)
(\div_{\pi_\tsigma}^\Omega V)(\om)\,\pi_\tsigma(d\om)\]
for all $F\in\FC$, provided it exists (i.e., provided
\[ F\mapsto\int_{\OM}\la V_\om,\nabla^\Omega
F(\om)\ra_{T_\om\big(\OM\big)}\pi_\tsigma(d\om)\]
is continuous on $L^2(\pi_\tsigma)$).
}\end{define}

A class of smooth vector fields on $\OM$ for which the divergence can be
computed in an explicit form is described in the following proposition.

\begin{prop}\label{prop4} For any vector field
\[ V_\om(x)=\sum_{j=1}^N F_j(\om)(v_j(x),u_j(x)),\qquad \om\in\OM,\ x\in X,\]
with $F_j\in\FC$\rom, $(v_j,u_j)\in\agot$\rom, $j=1,\dots,N$\rom, we have
\begin{align*}
(\div^\Omega_{\pi_\tsigma}V)(\om)&=\sum_{j=1}^N\big(\nabla_{(v_j,u_j)}^\Omega
F_j\big)(\om)+\sum_{j=1}^N B_{(v_j,u_j)}^{\pi_{\tsigma}}(\om)F_j(\om)\\
&=\sum_{j=1}^N
\la\nabla^\Omega F_j(\om), (v_j,u_j)\ra_{T_\om\big(\OM\big)}+\sum_{j=1}^N\la\beta^
\tsigma_{(v_j,u_j)}
,\om\ra F_j(\om).\end{align*}\end{prop}

\noindent{\it Proof}. Due to the linearity of $\nabla^\Omega$, it
is sufficient to consider the case $N=1$, i.e.,
$V_\om(x)=F_1(\om)(v(x),u(x))$. By Theorem~\ref{thIbP}, we have
for all $F_2\in\FC$
\begin{gather*}
-\int_{\OM}\la V_\om,\nabla^\Omega F_2
(\om)\ra_{T_\om\big(\OM\big)}\pi_{\tsigma}(d\om)=
-\int_{\OM}F_1(\om)\nabla_{(v,u)}^\Omega
F_2(\om)\,\pi_\tsigma(d\om)\\
=\int_{\OM}\big(\nabla^\Omega_{(v,u)}F_1\big)(\om)F_2(\om)\,\pi_\tsigma(d\om)+\int
_{\OM}F_1(\om)F_2(\om)B_{(v,u)}^{\pi_\tsigma}(\om)\,\pi_{\tsigma}(d\om),\end{gather*}
which yields
\begin{align*}
(\div_{\pi_\tsigma}^\Omega
V)(\om)&=\nabla^\Omega_{(v,u)}F_1(\om)+B^{\pi_\tsigma}_{(v,u)}(\om)F_1(\om)\\
&=\la\nabla^\Omega
F_1(\om),(v,u)\ra_{T_\om\big(\OM\big)}+\la\beta^\tsigma_{(v,u)},\om\ra
F_1(\om).\quad \blacksquare\end{align*}

\begin{rem}\label{rem1}\rom{
Extending the definition of $B^{\pi_\tsigma}$ in \eqref{16} to the class of
vector fields $V=\sum_{j=1}^N F_j\otimes (v_j,u_j)$ by
\[ B_V^{\pi_\tsigma}(\om):=\sum_{j=1}^N\la\beta^\tsigma_{(v_j,u_j)},\om\ra
F_j(\om)+\sum_{j=1}^N
\big(\nabla^\Omega_{(v_j,u_j)}F_j\big)(\om),\]
we obtain that
\[\div^\Omega_{\pi_\sigma}{\scriptstyle\bullet}=B^{\pi_\tsigma}_{\bullet}.\]
In particular, if $(v,u)\in\agot$, it follows,  for the ``constant'' vector field
$V_\om\equiv(v,u)$ on $\OM$, that
\[ \div^\Omega_{\pi_\tsigma}(v,u)(\om)=\la\div^{\XM}_\tsigma(v,u),\om\ra,\]
where $\div^\XM_\tsigma(v,u)=\beta^\tsigma_{(v,u)}$ is the
divergence on $\XM$ of $(v,u)$  w.r.t.\ $\tsigma$:
\begin{align*}
&\int_{\XM}\la\nabla^\XM\fii(x,m),(v(x),u(x))\ra_{T_{(x,m)}(\XM)}\,\tsigma(dx,dm)\\
&\qquad=-\int_\XM\fii(x,m)\big(\div_\tsigma^{\XM}(v,u)\big)(x,m)\,\tsigma(dx,dm),\qquad
\fii\in{\frak D}.\end{align*}
}\end{rem}

\subsection{Integration by parts characterization}

In the works \cite{AKR,AKR2} it was shown that the mixed Poisson measures are
exactly the ``volume elements'' corresponding to the differential geometry on the
configuration space $\Gamma_X$. Now, we wish to prove that an analogous
statement holds true in our case of $\OM$ for mixed marked Poisson measures.

We start with a lemma that describes $\tsigma$ as the unique (up
to a constant) measure on $\XM$ with respect to which the
divergence $\div^\XM_\tsigma$ is the dual operator of the gradient
$\nabla^{\XM}$.

\begin{lem}\label{lemcharacter}
Let the conditions \eqref{condition} and \eqref{neues} hold\rom.
Then\rom, for every $\Lambda\in\Oc(X)$ the measures
$z\tsigma$\rom, $z> 0$\rom, are the only positive Radon measures
$\xi$ on $\Lext$ such that $\div_\tsigma^\XM$ is the dual operator
on $L^2(\Lext;\xi)$ of $\nabla^\XM$ when considered with the
domains $V_0(\Lambda)\times C^\infty_0(\Lambda;\got)$\rom,
resp\rom.\ $C^\infty_{\text{\rom{0,b}}}(\Lext)$ \rom(i\rom.e\rom.,
the set of all $(v,u)\in\agot$\rom, resp\rom.\ $\varphi\in{\frak
D}$ with support in $\Lambda$\rom, resp\rom.\ $\Lext$\rom{).}
\end{lem}

\noindent {\it Proof}. In virtue of the conditions \eqref{condition} and
\eqref{neues},
the lemma is
obtained in complete analogy with Remark~4.1 (iii) in \cite{AKR2}. Indeed, let
$q_1(x,m)$ and $q_2(x,m)$ be two densities w.r.t.\ $\nu\otimes \lambda$ for which
the logarithmic derivatives coincide. Then, we get
\begin{align*}
\nabla_v^X\log q_1(x,m)&=\nabla_v^X\log q_2(x,m),\qquad v\in V_0(X),\\
\nablat_u^M\log q_1(x,m)&=\nablat_u^M\log q_2(x,m),\qquad u\in C_0^\infty(\Lambda;\got),
\,\text{$\nu\otimes\lambda$-a.s.},
\end{align*}
which yields respectively
\begin{align*}q_1(x,m)&=q_2(x,m)c(m),\\
q_1(x,m)&=q_2(x,m)\tilde c(x)\quad\text{$\nu\otimes\lambda$-a.s.}
\end{align*}
Therefore, $q_1(x,m)=\operatorname{const} q_2(x,m)$ $\nu\otimes\lambda$-a.s.\quad $\blacksquare$

Let $\kappa$ be a probability measure on $(\R_+,{\cal B}(\Rp))$.
Then, we define a mixed marked Poisson measure as follows:
\begin{equation}\label{salko} \mu_{\kappa,\tsigma}
=\int_\Rp\pi_{z\tsigma}\,\kappa(dz).\end{equation}
Here, $\pi_{0\tsigma}$ denotes the Dirac measure on $\OM$ with mass in $\om=\{\varnothing\}$.
 Let ${\cal M}_l(\OM)$, $l\in[1,\infty)$, denote the set of all probability
 measures on $(\OM,{\cal B}(\OM))$ such that
\[\int_{\OM}|\la f,\om\ra|^l\,\mu(d\om)<\infty\quad \text{for all }f\in
 C_{\text{0,b}}(\XM),\ f\ge0.\]
Clearly, $\mu_{\kappa,\tsigma}\in {\cal M}_l(\OM)$ if and only if
\begin{equation}\label{30}
\int_{\Rp}z^l\,\kappa(dz)<\infty.\end{equation}
We define $(\text{IbP})^\tsigma$ to be the set of all
$\mu\in {\cal
  M}_1(\OM)$ with the property that $\om\mapsto \la\beta^\tsigma_{(v,u)},\om\ra$
is $\mu$-integrable for all $(v,u)\in\agot$ and which satisfy \eqref{11} with $\mu$
replacing $\pi_\tsigma$ for all $F_1,F_2\in\FC$, $(v,a)\in\got$.
 We note that \eqref{11} makes sense only for such measures and that
 $B_{(v,u)}^{\pi_\tsigma}$ depends only on $\tsigma$ not on
 $\pi_\tsigma$. Obviously, since $\nabla^\XM_{(v,u)}$ obeys the product rule
 for all $(v,u)\in\agot$, we can always take $F_2\equiv
 1$. Furthermore, $(\text{IbP})^{\tsigma}$ is convex.

\begin{th}\label{hahaha} Let the condition \eqref{condition} and \eqref{neues} be satisfied\rom.
Then\rom, the
  following conditions are equivalent\rom:\\
\rom{(i)} $\mu\in({\mathrm IbP})^\tsigma$\rom;\\
\rom{(ii)} $\mu=\mu_{\kappa,\tsigma}$ for some probability measure $\kappa$ on
$(\Rp,{\cal B}(\Rp))$ satisfying \eqref{30} with $l=1$\rom. \end{th}

\noindent{\it Proof}. The part (ii)$\Rightarrow$(i) is trivial.
The proof of  (i)$\Rightarrow$(ii) goes along absolutely analogously to that
in the particular
case where $G=M=\R_+$, see \cite{KLU1}.
\vspace{2mm}\quad $\blacksquare$

As a direct consequence of Theorem~\ref{hahaha}, we obtain

\begin{cor} The extreme points of $(\operatorname{IbP})^\tsigma$ are
  exactly $\pi_{z\tsigma}$\rom, $z\ge0.$
\end{cor}

\subsection{A lifting of the geometry}

Just as in the case of the geometry on the configuration space, we can
present
an interpretation of the formulas obtained in subsections~3.1--3.3 via a
simple
``lifting rule.''

Suppose that $f\in C_{\text{0,b}}(\XM)$, or more generally $f$ is an arbitrary  measurable
function on $\XM$ for which there exists (depending on $f$) $\Lambda\in\Bc(X)$
such that $\supp f\subset\Lext$. Then, $f$  generates a (cylinder) function on $\OM$
by the formula
\[ L_f(\om):=\la f,\om\ra,\qquad \om\in\OM.\]
We will call $L_f$ the lifting of $f$.

As before, any vector field $(v,u)\in\agot$,
\[(v,u)\colon X\ni x\mapsto (v(x),u(x))\in T_{(x,m)}(\XM)=T_x(X)\dotplus
\got,\]
can be considered as a vector field on $\OM$ (the lifting of $(v,u)$),
which we
denote by $L_{(v,u)}$:
\[L_{(v,u)}\colon \OM\ni\om=\{\gamma,s\}\mapsto\{x\mapsto(v(x),u(x))\}\in
T_\om(\OM)=L^2(X\to T(X)\dotplus\got\,;\gamma).\]
For $(v_1,u_1), (v_2,u_2)\in\agot$, the formula \eqref{9.1}
can be written as follows:
\[\big\la L_{(v_1,u_1)},L_{(v_2,u_2)}\big\ra_{T_\om\big(\OM\big)}
=L_{\la(v_1,u_1),(v_2,u_2)\ra_{T(\XM)}}(\om),\]
i.e., the scalar product of lifted vector fields is computed as the
lifting of
the scalar product
\[\la(v_1(x),u_2(x)),(v_2(x),u_2(x))\ra_{T_{(x,s_x)}(\XM)}=f(x).\]
This rule can be used as a definition of the tangent space $T_\om\big(\OM\big)$.

The formula \eqref{8} has now the following interpretation:
\begin{equation}\label{20}
\big(\nabla^\Omega_{(v,u)}L_\fii\big)(\om)=L_{\nabla^\XM_{(v,u)}\fii}
(\om),\qquad
\fii\in{\frak D},\ \om\in\OM,
\end{equation}
and the ``lifting rule'' for the gradient is given by
\begin{equation}\label{17}(\nabla^\Omega L_\fii)(\gamma,s)\colon\gamma\ni
  x\mapsto\nabla^{\XM}\fii( x,s_x).\end{equation}

As follows from \eqref{16}, the logarithmic derivative
$B^{\pi_\tsigma}_{(v,u)}\colon\OM\to\R$
is obtained via the lifting procedure of the corresponding logarithmic
derivative $\beta^\tsigma_{(v,u)}\colon\XM\to\R$, namely,
\[ B^{\pi_\tsigma}_{(v,u)}(\om)=L_{\beta^\tsigma_{(v,u)}}(\om),\]
or equivalently, one has for the divergence of a lifted vector field:
\begin{equation}\label{18}\div^\Omega_{\pi_\tsigma}L_{(v,a)}=L_{\div^\XM_\tsigma}(v,a).\end{equation}

We underline that by \eqref{20} and \eqref{17} one recovers the action of
$\nabla^\Omega_{(v,a)}$ and $\nabla^\Omega$ on all functions from $\FC$
algebraically from requiring the product or the chain rule to hold. Also, the
action of $\div_{\pi_\tsigma}^\Omega$ on more general cylindrical vector fields
follows as in Remark~\ref{rem1} if one assumes the usual product rule for
$\div_{\pi_\tsigma^\Omega}$ to hold.

\section[Representations of the Lie algebra $\frak a$ of the group $\Agot$]{Representations
of the Lie algebra $\frak a$ of the group $\Agot$}

Using the $\Agot$-quasiinvariance of $\pi_{\tsigma}$, we can
define the unitary representation of the group
$\Agot=\Diff\underset{\alpha}{\times}G^X$ in the space
$L^2(\pi_{\tsigma})$. Namely, for $a\in\Agot$, we define the
unitary operator
\[\big(V_{\pi_\tsigma}(a)F\big)(\om):=F(a(\om))\sqrt{\frac{da^{-1*}\pi_\tsigma
    }{d\pi_\tsigma}(\om)},\qquad F\in L^2(\pi_\tsigma).\]
Then, we have
\[ V_{\pi_\tsigma}(a_1)V_{\pi_\tsigma}(a_2)=V_{\pi_\tsigma}(a_1a_2),\qquad
a_1,a_2\in\Agot.\]

As has been noted in Introduction, this representation is
reducible, cf.\ \cite{KLU1}

As in subsec.~3.1, to any vector field $v\in V_0(X)$ there corresponds a
one-parameter subgroup of diffeomorphisms $\psi_t^v$, $t\in\R$. It generates a
one-parameter unitary group
\[V_{\pi_\tsigma}(\psi_t^v):=\exp[itJ_{\pi_\tsigma}(v)],\qquad t\in\R,\]
where $J_{\pi_\tsigma}(v)$ denotes the selfadjoint generator of
this group. Analogously, to  a subgroup $\eta_t^u$, $u\in
C^\infty_0(X;\got)$, there corresponds a one-parameter unitary
group
\[ V_{\pi_\tsigma}(\eta_t^u):=\exp [it I _ {\pi_\tsigma}(u)]\]
with a generator $I_{\pi_\tsigma}(u)$.

\begin{prop}\label{prop5}
For any $v\in V_0(X)$ and $u\in C^\infty_0(X;\got)$\rom, the following operator equalities
on the domain $\FC$ hold\rom:
\begin{align*}
J_{\pi_\tsigma}(v)&=\frac 1i\,\nabla_v^\Omega+\frac{1}{2i}\,B_v^{\pi_\tsigma},\\
I_{\pi_\tsigma}(u)&=\frac1i\,\nabla_u^\Omega+\frac1{2i}\,B_u^{\pi_\tsigma}.
\end{align*}\end{prop}

\noindent{\it Proof}. These equalities follow immediately from the definition of
the directional derivatives $\nabla_v^\Omega$ and $\nabla_a^\Omega$,
Theorem~\ref{thIbP}, and the form of the operators $V_{\pi_\tsigma}(\psi_t^v)$
and $V_{\pi_\tsigma}(\theta_t^u)$.\quad $\blacksquare$ \vspace{2mm}

For any $(v,u)\in\agot$, define an operator
\[{\cal R}_{\pi_\tsigma}(v,u):= J_\pit(v)+I_\pit(u).\]
By Proposition~\ref{prop5},
\[{\cal R}_\pit(v,u)=\frac1i\nabla_{(v,u)}^\Om+\frac1{2i}B^\pit_{(v,u)}.\]
We wish to derive now a commutation relation between these operators.

\begin{lem}\label{lem3.3} The Lie-bracket $[(v_1,u_1),(v_2,u_2)]$ of the vector fields
  $(v_1,u_1)$\rom, $(v_2,u_2)\in \agot$\rom, i\rom.e\rom{.,} a
  vector field from $\agot$ such that
\[
  \nabla^{\XM}_{[(v_1,u_1),(v_2,u_2)]}=\nabla^{\XM}_{(v_1,u_1)}
\nabla^\XM_{(v_2,u_2)}-
\nabla^\XM_{(v_2,u_2)}\nabla^\XM_{(v_1,u_1)}\quad \text{\rom{on} }{\frak D},\]
is given by
\[[(v_1,u_1),(v_2,u_2)]=([v_1,v_2],\nabla^X_{v_1}u_2-\nabla^X_{v_2}u_1+ [u_1,u_2]),\]
where $[v_1,v_2]$ is the Lie-bracket of the vector fields $v_1$\rom, $v_2$ on
$X$\rom,
$$[u_1,u_2](x)=[u_1(x),u_2(x)]$$
 \rom(the latter being the Lie-bracket on
$\got$ of $u_1(x),u_2(x)\in\got$\rom{),} and $\nabla^X_v u$ is the derivative in direction $v$ of
a $\got$-valued function $u$ on $X$\rom.
\end{lem}

\noindent{\it Proof}. First, we have on $\frak D$:
\begin{equation}\label{raz}
\nabla^X_{v_1}\nabla^X_{v_2}-\nabla^X_{v_2}\nabla^X_{v_1}=
\nabla^X_{[v_1,v_2]}, \qquad v_1, v_2\in V_0(X).\end{equation}
Next, using \eqref{bum},  $$\nablat_u^M f(x,m)=\la \nabla^G\hat
f(x,e,m),u(x)\ra_\got,\qquad \hat f(x,g,m):=f(x,\theta(g,m)),$$
and so
\begin{align}
&(\nablat^M_{u_1}\nablat^M_{u_2}-\nablat^M_{u_2}\nablat^M_{u_1}
)f(x,m)\notag\\ &\qquad= \la\nabla^G\hat
f(x,e,m),[u_1(x),u_2(x)]\ra_\got \notag\\ &\qquad
=\nablat^M_{[u_1,u_2]}f(x,m),\qquad u_1,u_2\in C_0^\infty
(X;\got).\label{dva}
\end{align}

Finally,
\begin{gather}
(\nabla^X_v\nablat^M_u-\nablat_u^M\nabla^X_v)f(x,m)\notag\\
=\la\nabla^X\la\nabla^G\hat f(x,e,m),u(x)\ra_\got,v(x)\ra_{T_x(X)}
\notag\\ \mbox{}-\la\nabla^G\la\nabla^X\hat
f(x,e,m),v(x)\ra_{T_x(X)},u(x)\ra_\got \notag\\
=\la\nabla^X\nabla^G\hat f(x,e,m),v(x)\otimes
u(x)\ra_{T_x(X)\otimes \got}+\la\nabla^G\hat f(x,e,m),\nabla_v^X
u(x)\ra_\got\notag\\ \mbox{}-\la\nabla^G\nabla^X\hat
f(x,e,m),u(x)\otimes v(x)\ra_{\got \otimes T_x(X)}\notag\\
=\la\nabla^G\hat f(x,e,m),\nabla^X_v
u(x)\ra_\got=\nablat^M_{\nabla ^X_vu}f(x,m),\notag\\ v\in V_0(X),\
u\in C_0^\infty(X;\got). \label{tri}\end{gather} The equalities
\eqref{raz}--\eqref{tri} yield the lemma.\quad $\blacksquare$

\begin{prop}\label{prop3.3}
For arbitrary $(v_1,u_1),\, (v_2,u_2)\in\agot$\rom, the
following operator equality holds on $\FC$\rom:
$$
[{\cal R}_\pit(v_1,u_1),{\cal R}_\pit(v_2,u_2)]={\cal
  R}_\pit([(v_1,u_1),(v_2,u_2)]).$$
In particular\rom,
\begin{alignat*}{2}
[J_\pit(v_1),J_\pit(v_2)]&=-iJ_\pit([v_1,v_2]),&&\qquad v_1,v_2\in
V_0(X),\\ [I_\pit(u_1),I_\pit(u_2)]&=-I_\pit([u_1,u_2]) ,&&\qquad
u_1,u_2\in C_0^\infty(X;\got),\\
[J_\pit(v),I_\pit(u)]&=-iI_\pit(\nabla_v^Xu),&&\qquad v \in
V_0(X),\, u\in C_0^\infty(X;\got).\end{alignat*}
\end{prop}

\noindent {\it Proof}.
First we note that Lemma~\ref{lem3.3} and \eqref{8}
immediately imply
\[ \nabla^\Omega_{(v_1,u_1)}\nabla^\Omega_{(v_2,u_2)}-
\nabla^\Omega_{(v_2,u_2)}\nabla^\Om_{(v_1,u_1)}
=\nabla^\Om_{[(v_1,u_1),(v_2,u_2)]}\quad \text{on }\FC.\]
Therefore, by using the chain rule, we conclude that the lemma will be proved if
we show that
\begin{equation}\label{zyrk}
\nabla^\Om_{(v_1,u_1)}B^\pit_{(v_2,u_2)}-\nabla^\Om_{(v_2,u_2)}
B^\pit_{(v_1,u_1)}=B^\pit_{[(v_1,u_1),(v_2,
  u_2)]}\quad \text{$\pit$-a.e.}\end{equation}
But upon the representation
\[ B^\pit_{(v,u)}(\om)=\la \nabla^X_v\log q+\nablat^M_u\log q+
\div^X v+ \la\nabla^G p^\lambda(e,m),-u(x)\ra_\got,\om\ra\] and
Remark~\ref{kiii}, we easily derive \eqref{zyrk} again from
Lemma~\ref{lem3.3}.\quad $\blacksquare$\vspace{2mm}

Thus, the operators ${\cal R}_{\pi_\tsigma}(v,u)$, $(v,u)\in\agot$,
give a marked
Poisson space representation of the Lie algebra $\agot$ of the group
$\Agot$.

\section[Intrinsic Dirichlet forms on marked Poisson spaces]{Intrinsic Dirichlet forms on marked Poisson spaces}

\subsection{Definition of the intrinsic Dirichlet form}

From now on, the underlying space of ``nice functions'' on $\XM$
will be instead of $\frak D$ the space ${\frak
D}_0:=C_0^\infty(\XM)$ consisting of all $C^\infty$ functions with
compact support in $\XM$. Evidently, ${\frak D}_0$ is a subset  of
$\frak D$ and in the case where $M$ is itself compact ${\frak
D}_0=\frak D$. Absolutely analogously to $\FC$ one constructs the
set $\FCO$($\subset\FC$), which is dense in $L^2(\pi_\tsigma)$. By
$\FPO$ we denote the set of all cylinder functions of the form
\eqref{5} in which the functions $\fii_1,\dots,\fii_N$ belong to
${\frak D}_0$ and the generating function $g_F$ is a polynomial on
$\R^N$, i.e., $g_F\in{\cal P}(\R^N)$. Finally, in the same way we
introduce $\FCpO$ where $g_F\in C_{\mathrm p}^\infty(\R^N)$ (:=the
set of all $C^\infty$-functions $f$ on $\R^N$ such that $f$ and
its partial derivatives of any order are polynomially bounded).

We have obviously
\begin{align*}
\FCO&\subset\FCpO,\\
\FPO&\subset\FCpO,\end{align*}
and these are algebras with respect to the usual operations. The existence of
the Laplace transform $\ell_{\pi_\tsigma}(f)$ for each
$f\in C_{\mathrm
  0}(\XM)$ implies, in particular, that $\FCpO\subset L^2(\pi_\tsigma)$.

\begin{define}\rom{
For $F_1,F_2\in\FCpO$, we introduce a pre-Dirichlet form as
\begin{equation}\label{4.1}
\EOM(F_1,F_2)=\int_{\OM}\la\nabla^\Om F_1(\om),\nabla^\Om
F_2(\om)\ra_{T_\om(\OM)}\,\pi_{\tsigma}(d\om).\end{equation}
}\end{define}

Note that, for all $F\in\FCpO$, the formula \eqref{gradient} is still valid and
therefore, for $F_1=g_{F_1}(\la\fii_1,\cdot\ra,\dots,\la\fii_N,\cdot\ra)$ and
$F_2=g_{F_2}(\la\xi_1,\cdot\ra,\dots,\la\xi_K, \cdot\ra)$ from $\FCpO$, we have
\begin{gather}
\la\nabla^\Om F_1(\om),\nabla^\Om F_2(\om)\ra_{T_\om(\OM)}=\notag\\
=\sum_{j=1}^N\sum_{k=1}^K\frac{\di g_{F_1}}{\di
  r_j}(\la\fii_1,\om\ra,\dots,\la\fii_N,\om\ra)\frac{\di g_{F_2}}{\di
  r_k}(\la\xi_1,\om\ra,\dots,\la\xi_K,\om\ra)\times\notag\\
\times\int_X\la\nabla^{\XM}
\fii_j(x,s_x),\nabla^{\XM}\xi_k(x,s_x)\ra_{T_{(x,s_x)}(\XM)}\,\gamma(dx)\notag\\
=\sum_{j=1}^N\sum_{k=1}^K\frac{\di g_{F_1}}{\di
  r_j}(\la\fii_1,\om\ra,\dots,\la\fii_N,\om\ra)\frac{\di g_{F_2}}{\di
  r_k}(\la\xi_1,\om\ra,\dots,\la\xi_K,\om\ra)\times\notag\\
\times\la\la\nabla^{\XM}\fii_j,\nabla^{\XM}\xi_k\ra_{T(\XM)},\om\ra.\label{sex}
\end{gather}
Since for $\fii,\xi\in{\frak D}_0$, the function
\begin{gather*}\la\nabla^\XM\fii(x,m),\nabla^\XM\xi(x,m)
\ra_{T_{(x,m)}(\XM)}=\\
=\la\nabla^X\fii(x,m),\nabla^X\xi(x,m)\ra_{T_x(X)}
+\la\nablat^{M}\fii(x,m),\nablat^M\xi(x,m)\ra_\got\end{gather*}
belongs to ${\frak D}_0$, we conclude that
\[\la\nabla^\Om F_1(\cdot),\nabla^\Om(\cdot)F_2(\cdot)\ra_{T(\OM)}\in
L^1(\pi_\tsigma),\qquad F_1,F_2\in\FCpO,\]
and so \eqref{4.1} is well defined.

We will call $\EOM$ the intrinsic pre-Dirichlet form corresponding to the marked
Poisson measure $\pi_\tsigma$ on $\OM$. In the next subsection we will prove the
closability of $\EOM$.

\subsection{Intrinsic Dirichlet operators}

We start with introducing the pre-Dirichlet operator corresponding to
the measure $\tsigma$ on $\XM$ and to the gradient $\nabla^\XM$:
\begin{equation}
{\cal
  E}^\XM_\tsigma(\fii,\xi):=\int_\XM\la\nabla
^\XM\fii(x,m),\nabla^\XM\xi(x,m)\ra_{T_{(x,m)}(\XM)}\,\tsigma (dx,
dm),\label{4.5}\end{equation} where $\fii,\,\xi\in{\frak D}_0$.
This form is associated with the Dirichlet
  operator
\begin{equation}\label{z4.4}H_\tsigma^\XM:=H_\tsigma^X+H_\tsigma^M\end{equation}
 on ${\frak D}_0$ which
  satisfies
\begin{equation}\label{z4.5} {\cal
  E}^\XM_\tsigma(\fii,\xi)=(H^\XM_\tsigma\fii,\xi)
_{L^2(\tsigma)},\qquad
  \fii,\,\xi\in{\frak D}_0.\end{equation}
Here, $H^X_\tsigma$ and $H^M_\tsigma$ are the Dirichlet operators of
$\nabla^X$ and $\nablat^M$, respectively. Evidently,
\begin{equation}
H^X_\tsigma\fii(x,m)=-\Delta^X\fii(x,m)-\la\nabla^X\log
q(x,m),\nabla^X\fii(x,m)\ra_{T_x(X)},\label{4.6}\end{equation}
where $\Delta^X$  denotes the Laplace-Beltrami operator corresponding
to $\nabla^X$.

\newcommand{\divtM}{\widetilde{\div}{}^M}

Let us calculate the operator $H_\tsigma^M$. Suppose $f\in{\frak
D}_0$ and $W\in C_0(\XM;\got)$. Analogously to
Remark~\ref{brandnewremark}, we conclude
\begin{equation}\label{shudfzrwe} \la\nablat^M f(x,m),W(x,m)\ra_\got
=\la\nabla^M f(x,m),(RW)(x,m)\ra_{T_m(M)}, \end{equation} where
$RW\in C_0^\infty(\XM;TM)$ is given by
\begin{equation}
\label{urgan} X\times M\ni (x,m)\mapsto(RW)(x,m):=\frac
d{dt}\theta(\exp(tW(x,m)),m)\big|_{t=0}\in T_mM.\end{equation}
 Therefore, using the integration
by parts formula on $M$ for a vector field with a compact support,
we get
\begin{gather*}
\int_\XM \la\nablat^M f(x,m),W(x,m)\ra_\got\,\tsigma(dx,dm)\\
=-\int_{\XM} f(x,m)\big[\div^M(RW)(x,m)\\ \text{}+\la\nabla^M\log
q(x,m), (RW)(x,m)\ra_{T_m(M)}\big]\,\tsigma(dx,dm)\\
=-\int_{\XM}f(x,m)\big[\div^M(RW)(x,m)+\la\nablat^M\log q(x,m),
W(x,m)\ra_\got\big]\,\tsigma(dx,dm),\end{gather*} where $\div^M$
is the divergence on $M$ with respect to the usual gradient
$\nabla^M$ and the measure $\lambda$. Thus, the divergence
$\widetilde{\div}{}^M_{\tsigma}$ on $\XM$ w.r.t.\ the gradient
$\widetilde\nabla{}^M$ and the measure $\tsigma$ is given by
$$\widetilde{\div}{}^M_{\tsigma}W(x,m)=\div^M
(RW)(x,m)+\la\nablat^M\log q(x,m), W(x,m)\ra_\got.$$ In
particular, the divergence $\widetilde{\div}{}^M$  w.r.t.\ the
 measure
$\nu(dx)\,\lambda(dm)$ equals
\begin{equation}\widetilde{\div}{}^M W(x,m)=\div^M
(RW)(x,m).\label{sprth}\end{equation}

 It is easy to see that, for $f\in{\frak D}_0$, $W=\nablat^M f\in
C_0^\infty( \XM;\got)$, and so we have finally
\begin{equation}
\label{4.6prime} H_\tsigma^M f=
\widetilde\div{}^M\widetilde\nabla{}^Mf= -\tilde \Delta^M
f-\la\nablat^M\log q,\nablat^M f\ra_\got, \qquad f\in{\frak
D}_0,\end{equation} where
\begin{equation}
\label{sezam} \tilde \Delta^M f=\divtM\nablat^M
f:=\div^M(R(\nablat^Mf)).
\end{equation}

The closure of the form ${\cal E}^\XM_\tsigma$ on
$$L^2(\XM;\tsigma)=:L^2(\tsigma)$$ is denoted by $({\cal
E}^\XM_\tsigma, D({\cal E}^\XM_\tsigma))$. This form generates a
positive selfadjoint operator in $L^2(\tsigma)$ (the so-called
Friedrichs extension of $H_\tsigma^\XM$, see e.g.\ \cite{BKR}).
For this extension we preserve the notation $\HXM$ and denote the
domain by $D(\HXM)$.

Let us introduce a differential operator $\HOM$ on the domain $\FCO$ which is
given on any $F\in\FCO$ of the form \eqref{5} by the formula
\begin{align}
(\HOM F)(\om):&=-\sum_{j,k=1}^N\frac{\di^2 F}{\di r_j\di r_k
  }(\la\fii_1,\om\ra,\dots,\la\fii_n,\om\ra)
\la\la\nabla^\XM\fii_j,\nabla^\XM\fii_k\ra_{T(\XM)},\om\ra \notag\\
&\quad+\sum_{j=1}^N\frac{\di F}{\di r_j}(\la\fii_1,\om\ra,\dots,\la\fii_n,
\om\ra)\la\HXM\fii_j,\om\ra.
\label{4.7}\end{align}
Since
\[ \la\nabla^\XM\log q,\nabla^\XM\fii_j\ra_{T(\XM)}\in L^2(\tsigma)\cap
L^1(\tsigma)\]
(see condition \eqref{condition}), the r.h.s.\
of \eqref{4.7} is well defined as an element of $L^2(\pi_\tsigma)$.
The  following theorem  implies, in particular, that
$\HOM$ is well defined as a linear operator on $\FCO$,
i.e., independently of the representation of $F$ as in \eqref{5}.

\begin{th} \label{th4.1}The operator $\HOM$ is associated with the intrinsic Dirichlet form
  $\EOM$ in the sense that\rom,  for all $F_1,F_2\in\FCO$
\begin{equation}\label{fr}\EOM(F_1,F_2)=(\HOM F_1,F_2)
_{L^2(\pi_\tsigma)},\end{equation}
or
\[ \HOM=-\div^\Om_{\pi_\tsigma}\nabla^\Om\quad {\mathrm on}\ \FCO.\]
We call $\HOM$ the intrinsic Dirichlet operator of the measure
$\pi_\tsigma$\rom. \end{th}

\begin{lem}\label{newlemma} For any $\fii\in{\frak D}_0$
and $W\in C_0^\infty(\XM;\got)$\rom, we have
$$\divtM(\fii W)(x,m)=\la\nablat^M\fii(x,m),W(x,m)\ra_\got+
\fii(x,m)\divtM W(x,m).$$
\end{lem}

\noindent{\it Proof}. By \eqref{shudfzrwe}, \eqref{urgan}, and
\eqref{sprth}
\begin{gather*}
\divtM(\fii W)(x,m)=\div ^M\Big[\frac
d{dt}\,\theta(\exp(t\fii(x,m)W(x,m)),m) \big|_{t=0}\Big]\\
=\div^M\Big[\fii(x,m)\,\frac
d{dt}\,\theta(\exp(tW(x,m)),m)\big|_{t=0}\Big]\\
=\la\nabla^M\fii(x,m),\frac
d{dt}\,\theta(\exp(tW(x,m)),m)\big|_{t=0} \ra_{T_m(M)}\\
\mbox{}+\fii(x,m)\div^M\Big[\frac
d{dt}\theta(\exp(tW(x,m)),m)\big|_ {t=0}\Big]\\ =\frac
d{dt}\,\fii(x,\theta(\exp(tW(x,m)),m))\big|_{t=0}+\fii(x,m) \divtM
W(x,m)\\ =\la\nablat^M\fii(x,m),W(x,m)\ra_\got+\fii(x,m)\divtM
W(x,m).\quad\blacksquare
\end{gather*}
\vspace{1mm}

\noindent {\it Proof of Theorem\/} \ref{th4.1}. For  shortness of
notations we will prove the formula \eqref{fr} in the case where
$F_1,F_2\in\FCO$ are of the form
\[ F_1=g_{F_1}(\la\fii,\om\ra),\quad F_2=g_{F_2}(\la\xi,\om\ra).\]
However, it is a trivial step to generalize the proof to general
$F_1,F_2$.

Let $\Lambda\in\Oc(X)$ be chosen so that the supports of the
functions $\fii$ and $\xi$ are in $\Lext$. Then, by \eqref{4.1},
\eqref{sex}, and the construction of the marked Poisson measure
\begin{gather*}
\EOM(F_1,F_2)=\int_{\OM}g_{F_1}'(\la\fii,\om\ra)g'_{F_2}(\la\xi,\om\ra)
\la\la\nabla^\XM\fii,\nabla^\XM\xi\ra_{T(\XM)},\om\ra \,\pi_\tsigma(d\om)\\
=-e^{\tsigma(\Lext)}\sum_{n=1}^\infty \frac 1{n!}
\int_{\Lambda^n_{\mathrm
    \ext}}g'_{F_1}(\fii(x_1,m_1)+\cdots+\fii(x_n,
m_n))\\
\times
g_{F_2}'(\xi(x_1,m_1)+\cdots+\xi(x_n,m_n))\\
\times\bigg[\sum_{i=1}^n\la\nabla^{\XM}\fii(x_i,m_i),\nabla^\XM
\xi(x_i,m_i)\ra_{T_{(x_i,m_i)}(\XM)}\bigg]
\tsigma(dx_1,dm_1)\dotsm\tsigma(dx_1,dm_1)\\
=e^{-\tsigma(\Lext)}\sum_{n=1}^\infty \frac 1{n!}\int_{\Lambda^n_{\mathrm
    \ext}}\sum_{i=1}^n\la\nabla_i^\XM g_{F_1}(\fii(x_1,m_1)+\cdots+
\fii(x_n,m_n)),\\
\nabla_i^\XM
g_{F_2}(\xi(x_1,m_1)+\dots+\xi(x_n,m_n))\ra_{
  T_{(x_i,m_i)}(\XM)}\,\tsigma(dx_1,dm_1)\dotsm\tsigma(dx_n,dm_n),
\end{gather*}
where $\nabla_i^\XM$ denotes the $\nabla^\XM$ gradient in the
$(x_i,m_i)$ variables. Therefore, by using \eqref{4.6prime} and
Lemma~\ref{newlemma}, we proceed in the calculation of
$\EOM(F_1,F_2)$ as follows: \allowdisplaybreaks
\begin{gather*}
 =e^{-\tsigma(\Lext)}\sum_{n=1}^\infty\frac
1{n!}\int_{\Lambda^n_{\mathrm \ext}}\bigg[\sum_{i=1}^n
H_\tsigma^{(\XM)_i}g_{F_1}(\fii(x_1,m_1)+\cdots+\fii(x_n,m_n))\bigg]
\times\\ \times
g_{F_1}(\xi(x_1,m_1)+\cdots+\xi(x_n,m_n))\,\tsigma(x_1,m_1)
\dotsm\tsigma(dx_n,dm_n)\\ =-e^{\tsigma(\Lext)}\sum_{n=1}^\infty
\frac 1{n!}\int_{\Lambda^n_{\mathrm
    \ext}}\bigg[ \sum_{i=1}^n g_{F_1}''(\fii(x_1,m_1)+\cdots+\fii(x_n,m_n))
\times \\ \times \la
\nabla^\XM\fii(x_i,m_i),\nabla^\XM\fii(x_i,m_i)\ra_{T_{(x_i,m_i)}(\XM)}\\
\mbox{}+g'_{F_1}(\fii(x_1,m_1)+\cdots+\fii(x_n,m_n))
H_\tsigma^{\XM}\fii(x_i,m_i) \bigg]\times\\ \times
g_{F_2}(\xi(x_1,m_1)+\cdots+\xi(x_n,m_n))\,\tsigma(dx_1,dm_1)
\dotsm\tsigma(dx_n,dm_n)\\ =\int_{\OM}\HOM
F_1(\om)F_2(\om)\,\pi_\tsigma(d\om).\quad\blacksquare
\end{gather*}

\begin{rem}\label{rem4.2}\rom{
The operator $\HOM$ can be naturally extended to cylinder functions of the form
\[ F(\om):=e^{\la\fii,\om\ra},\qquad \fii\in{\frak D}_0,\, \om\in\OM,\]
since such $F$ belong to $L^2(\pi_\tsigma)$. We then have
\begin{equation}\label{4.9}\HOM
  e^{\la\fii,\om\ra}=\la\HXM\fii-|\nabla^\XM\fii|^2_{T(\XM)},\om\ra
 \, e^{\la\fii,\om\ra}.\end{equation}
}\end{rem}

As an immediate  consequence of Theorem~\ref{th4.1} we obtain

\begin{cor}\label{cor4.1}
$(\EOM,\FCO)$ is closable on $L^2(\pi_\tsigma)$\rom. Its closure
$(\EOM,D(\EOM))$ is associated with a positive definite selfadjoint
operator\rom, the Friedrichs extension of $\HOM$\rom, which we also denote by
$\HOM$ \rom(and its domain by $D(\HOM)$\rom{).}\end{cor}

Clearly, $\nabla^\Om$
 also extends to $D(\EOM)$. We denote this extension by $\nabla^\Om$.

\begin{cor}\label{cor4.2} Let
\begin{equation}\label{4.10}
\begin{gathered}
F(\om):=g_F(\la\fii_1,\om\ra,\dots,\la\fii_N,\om\ra),\qquad \om\in\OM,\\
\fii_1,\dots,\fii_N\in D({\cal E}^\XM_\tsigma),\ g_F\in C^\infty_{\mathrm
  b}(\R^N).\end{gathered}\end{equation}
Then $F\in D(\EOM)$ and
\[ (\nabla^\Om F)(\omega;x)=\sum_{j=1}^N
\frac{\di g_F}{\di
  r_j}(\la\fii_1,\om\ra,\dots,\la\fii_N,\om\ra)\nabla^\XM\fii_j(x,s_x).\]
\end{cor}

\noindent{\it Proof}. By approximation this is an immediate  consequence of
\eqref{gradient} and the fact that, for all $1\le i\le N$,
\begin{equation}\label{4.11}
\int\la |\nabla^\XM\fii_i|^2_{T(\XM)},\om\ra\,\pi_\tsigma(d\om)={\cal
  E}^\XM_\tsigma(\fii_i,\fii_i).\end{equation}
\begin{rem}\label{rem4.3}
\rom{
Let $\mu_{\nu,\tsigma}\in{\cal M}_2(\OM)$ be given as in \eqref{salko}. Then, by
Theorem~3.2, (ii)$\Rightarrow$(i), all results above are
valid with $\mu_{\nu,\tsigma}$ replacing $\pi_\tsigma$. By \eqref{4.7} we have
\[ \HOM=H^\Om_{\mu_{\nu,\tsigma}}\quad\text{on }\FCO.\]
We note that the r.h.s.\ of \eqref{4.7} only depends on $\tsigma$
and the Riemannian structure of $\XM$. The respective Friedrichs
extension on $L^2(\mu_{\nu,\tsigma})$ is again denoted by
$H^\Om_{\mu_{\nu,\tsigma}}$, however it does necessarily  not
coincide with $\HOM$. }\end{rem}

\subsection{The heat semigroup and ergodicity}

The results of this subsection are obtained absolutely analogously to the
corresponding results of the paper \cite{AKR}, so we omit the proofs.

\newcommand{\munu}{\mu_{\kappa,\tsigma}}
\newcommand{\HH}{H^\XM_{\tsigma}}
\newcommand{\TTT}{T^\Om_{\munu}(t)}
\newcommand{\HHH}{H^\Om_{\munu}}
\newcommand{\EE}{{\cal E}_\tsigma^\XM}
\newcommand{\EEE}{{\cal E}^\Om_{\munu}}

For $\munu\in {\cal M}_2(\OM)$ let $\TTT:=\exp(-t\HHH)$, $t>0$. Define
\[ E({\frak
  D}_1,\OM)=\operatorname{l.h.}\big\{\,\exp(\la\log(1+\fii),\cdot\ra)\mid\fii\in{\frak
  D}_1\,\big\},\]
where l.h.\ means the linear hull and
\begin{align*}
{\frak D}_1:=\big\{\,\fii\in D(\HH)&\cap L^1(\tsigma)\mid \HH\fii\in
L^1(\tsigma)\\
&\text{and }-\delta\le\fii\le0\text{ for some }\delta\in(0,1)\,\big\}
.\end{align*}

\begin{prop}\label{prop4.1}
Let $\munu$ be as in \eqref{salko}\rom. Assume that $\HH$ is conservative\rom,
i\rom.e\rom{.,}
\[\int_{\XM}(\HH\fii)(x,m)\,\tsigma(dx,dm)=0\]
for all $\fii\in D(\HH)\cap L^1(\tsigma)$ such that $\HH\fii\in
L^1(\tsigma)$\rom, and suppose
\linebreak
that $(\HH,{\frak D}_0)$ is essentially selfadjoint
on $L^2(\tsigma)$\rom. Then
\begin{equation}\label{peter}
\TTT\exp(\la\log(1+\fii),\cdot\ra)=\exp(\la\log(1+e^{-t\HH}\fii),\cdot\ra),\qquad
\fii\in{\frak D}_1,\end{equation}
$E({\frak D}_1,\OM)\subset D(\HHH)$\rom, and
\begin{align*}
&\HHH\exp(\la\log(1+\fii),\cdot\ra)\\
&\qquad =\la (1+\fii)^{-1}\HH\fii,\cdot\ra\exp(\la\log(1+\fii),\cdot\ra),\qquad
\fii\in{\frak D}_1.\end{align*}
\end{prop}

\begin{rem}\rom{
(i) The condition of essential selfadjointness of $\HH$ on ${\frak
D}_0$ is fulfilled if $X$ is complete and
$|\beta^\tsigma|_{T(\XM)}\in L^p_{\mathrm loc}(\XM;m\otimes
\lambda)$ for some $p\ge\operatorname{dim}(X)+1$.

(ii) Since $(\exp(-t \HH))_{t>0}$ is sub-Markovian (i.e.,
$0\le\exp(-t\HH)\fii\le1$ for all $t>0$ and $\fii\in L^2(\tsigma)$,
$0\le\fii\le1$), because $(\EE,D(\EE))$ is a Dirichlet form, by a simple
approximation argument Proposition~\ref{prop4.1}  implies that the equality
\eqref{peter} holds for $t>0$ and
all $\fii\in L^1(\tsigma)$, $-1<\fii\le0$.
}\end{rem}

\begin{th}\label{th4.2}
Let the conditions of Proposition~\rom{\ref{prop4.1}} hold\rom. Then $E({\frak
  D}_1,\OM)$ is an operator core for the Friedrichs extension $\HHH$ on
  $L^2(\munu)$\rom. \rom(In other words\rom:\linebreak $(\HHH,E({\frak D}_1,\OM))$ is
  essentially selfadjoint on $L^2(\munu)$\rom{.)}
\end{th}

\begin{th}\label{th4.3} Suppose that the conditions of
  Theorem~\rom{\ref{hahaha}}  and Proposition~\rom{\ref{prop4.1}} hold\rom. Then the
  following assertions are equivalent\rom:

\rom{(i)} $\munu=\pi_{z\tsigma}$ for some $z>0$\rom.

\rom{(ii)} $(\EEE,D(\EEE))$ is irreducible \rom(i\rom.e\rom{.,} for $F\in
D(\EEE)$\rom, $\EEE(F,F)=0$ implies that $F=\operatorname{const}$\rom{).}

\rom{(iii)} $(\TTT)_{t>0}$ is irreducible \rom(i\rom.e\rom{.,} if $G\in
L^2(\munu)$ such that $\TTT(GF)=G\TTT F$ for all $F\in L^\infty(\munu)$\rom,
$t>0$\rom, then $G=\operatorname{const}$\rom{).}

\rom{(iv)} If $F\in L^2(\munu)$ such that $\TTT F=F$ for
all $T>0$\rom, then $F=\operatorname{const}.$

\rom{(v)} $\TTT\not\equiv{\bf 1}$ and ergodic \rom(i\rom.e\rom{.,}
\[\int\bigg(\TTT F-\int F\,d\munu\bigg)^2d\munu\to0\quad \mathrm{as}\ t\to0\]
for all $F\in L^2(\munu)$\rom{).}

\rom{(vi)} If $F\in D(\HHH)$ with $\HHH=0$\rom, then $F=\operatorname{const}.$
\end{th}

\begin{rem}\rom{
Let us consider the  diffusion process $P$ on $\XM$ associated to
the Dirichlet form $(\EE,D(\EE))$. This process can be interpreted
as distorted Brownian motion on the manifold $\XM$. More
precisely, the diffusion of points $x\in X$ is associated to the
Dirichlet form of the measure $\sigma$,
 so that it is distorted Brownian motion on
$X$, and the diffusion of marks $s_x$, $x\in X$, is associated to the
$\nablat^M$-Dirichlet form of the measure $p(x,dm)$ on $M$.

The existence of a diffusion process ${\bf P}$ corresponding to
the Dirichlet form \linebreak $(\EEE,D(\EEE))$ follows from
\cite{MaR}, and its identification with the independent infinite
particle process (on $\XM$) may be proved by the same arguments as
in \cite{AKR}. By analogy with the case of the process $P$ on
$\XM$, one can call $\bf P$ distorted Brownian motion on $\OM$.
}\end{rem}

\section{Intrinsic Dirichlet operator and second quantization}

\renewcommand{\TTT}{T^\Om_{\pi_\tsigma}(t)}
\renewcommand{\HHH}{H^\Om_{\pi_\tsigma}}
\renewcommand{\EEE}{{\cal E}^\Om_{\pi_\tsigma}}
\newcommand{\nMP}{\nabla^{\mathrm MP}}
\newcommand{\ot}{\otimes}
\newcommand{\hot}{\widehat\otimes}

In this section, we want to describe the Fock space realization of the marked
Poisson spaces and show that $\HHH$ is the second quantization of the operator
$\HH$.

\subsection{Marked Poisson gradient and chaos decomposition}

Let us define another ``gradient'' on functions $F\colon \OM\to\R$, which has
specific useful properties on the marked Poisson space.

\begin{define}\rom{
For any $F\in\FCpO$ we define the marked Poisson gradient
$\nMP$ as
\[ (\nMP F)(\om,(x,m)):=F(\om+\eps_{(x,m)})-F(\om),\qquad \om\in\OM,\ (x,m)
\in \XM.\]
}\end{define}

Let us mention that the operation
\[\OM\ni\om\mapsto\om+\eps_{(x,m)}\in\OM\]
is a $\pi_\tsigma$-a.e.\ well-defined map because of the property
\[ \pi_\tsigma\big(
\{\om=(\gamma,s)\in\OM\mid x\in\gamma\}
\big)=0\]
for an arbitrary $x\in X$ (which easily follows from the construction of
$\pi_\tsigma$). We consider $\nMP$ as a mapping
\[\nMP\colon \FCpO\ni F\mapsto \nMP F\in L^2(\tsigma)\otimes L^2(\pi_\tsigma)\]
that corresponds to using the Hilbert space $L^2(\tsigma)$ as a tangent space at
any point $\om\in\OM$. Thus, for any $\fii\in{\frak D}_0$,
we can introduce the
directional derivative
\begin{align*}
(\nabla^{\mathrm MP}_\fii F)(\om)&=\la\nMP F(\om),\fii\ra_{L^2(\tsigma)}\\
&=\int_{\XM}(F(\om+\eps_{(x,m)})-F(\om))\fii(x,m)\,\tsigma(dx,dm)
.\end{align*}

The most important feature of the marked Poisson gradient is that
it produces (via a corresponding ``integration by parts formula'')
the orthogonal system of Charlier polynomials on $(\OM,{\cal
B}(\OM),\pi_\tsigma)$. Below, we describe this construction in
detail using the isomorphism between $L^2(\pi_\tsigma)$ and the
symmetric Fock space (see \cite{ItKu,KSS2,LRS})

Let ${\cal F}(L^2(\tsigma))$ denote the symmetric Fock space over
$L^2(\tsigma)$:
\[ {\cal F}(L^2(\tsigma)):=\bigoplus_{n=0}^\infty {\cal F}_n(L^2(\tsigma))n!,\]
where \begin{gather*} {\cal
F}_n(L^2(\tsigma)):=(L^2(\tsigma))^{\hot n}=\hat
L^2((\XM)^n,\tsigma^{\ot n}),\qquad n\in\N,\\ {\cal
F}_0(L^2(\tsigma)):=\R,\end{gather*} $\hot$ denoting the symmetric
tensor product.  Thus, for each $F=(f^{(n)})_{n=0}^\infty\in {\cal
F}(L^2(\tsigma))$
\[ \| F\|^2_{{\cal F}(L^2(\tsigma))}=\sum_{n=0}^\infty |f^{(n)}|_{\hat
  L^2(\tsigma^{\ot n})}n!.\]

By ${\cal F}_{\mathrm fin}({\frak D}_0)$ we denote the dense subset of ${\cal
  F}(L^2(\tsigma))$ consisting of finite sequences $(f^{(n)})_{n=0}^N$,
  $n\in\Z_+$, such that each $f^{(n)}$ belongs to
${\cal F}_n({\frak D}_0):={\mathrm
  a}.{\frak D}_0^{\hot n}$, the $n$-th symmetric algebraic tensor power of ${\frak
  D}_0$:
\[ {\mathrm a}.{\frak D}_0^{\hot n}:=\operatorname{l.h.}\{\fii_1\hot\dotsm\hot
  \fii_n\mid\fii_i\in{\frak D}_0\}.\]
In virtue of the polarization identity, the latter  set is spanned just by the
vectors of the form $\fii^{\ot n}$ with $\fii\in{\frak D}_0$.

Now, we define a linear mapping
\begin{equation}{\cal F}_{\mathrm fin}({\frak D}_0)\ni F=(f^{(n)})_{n=0}^N\mapsto
IF=(IF)(\om)=\sum_{n=0}^N Q_n (f^{(n)};\om)\in\FPO
\label{mapping}\end{equation}
by using the following recursion relation:
\begin{align}
Q_{n+1}(\fii^{\ot(n+1)};\om)&=Q_n(\fii^{\ot
  n};\om)(\la\om,\fii\ra-\la\fii\ra_{\tsigma})\notag\\ &\quad-
nQ_n(\fii^{\ot(n-1)}\hot(\fii^2),\om)-nQ_{n-1}(\fii^{\ot(n-1)};\om)\la\fii^2\ra_{\tsigma},\notag\\
&\quad Q_0(1,\om)=1,\quad \fii\in{\frak D}_0.\label{ser}
\end{align}
Here, we have set $\la\fii\ra_{\tsigma}:=\int\fii\,d\tsigma$.
Notice that, since ${\frak
  D}_0$ is an algebra under pointwise multiplication of functions, the latter
  definition is correct.

  It is not hard to see that the mapping \eqref{mapping} is
  one-to-one. Moreover, the following proposition holds:

\begin{prop}\label{prop5.1} The mapping \eqref{mapping} can be extended by
  continuity to a unitary isomorphism between the spaces ${\cal
  F}(L^2(\tsigma))$ and $L^2(\pi_\tsigma).$\end{prop}

For each $\fii\in{\frak D}_0$, let us define the creation and annihilation
operators in ${\cal F}(L^2(\tsigma))$ by
\[ a^+(\fii)\psi^{\ot n}=\fii\hot\psi^{\ot n},\qquad a^-(\fii)\psi^{\ot
  n}=n(\fii,\psi)_{L^2(\tsigma )}\psi^{\ot (n-1)},\qquad \psi\in{\frak D}_0.\]

We will denote by the same letters the images of these operators under the
unitary $I$.

\begin{prop}\label{prop5.2} We have\rom, for each $\fii\in{\frak D}_0,$
\[ a^-(\fii)=\nabla^{\mathrm MP}_\fii,\qquad a^+(\fii)=\nabla_\fii^{{\mathrm
    MP}*}.\]
In particular\rom,
\[ Q_n(\fii_1\hot\dotsm\hot\fii_n;\om)=(\nabla^{{\mathrm MP}*}_{\fii_1}\dotsm
\nabla^{{\mathrm MP}*}_{\fii_n}{\bf 1})(\om),\qquad \om\in\OM.\]
\end{prop}

\newcommand{\Exp}{\operatorname{Exp}}

Finally, for each $\fii\in{\frak D}_0$ we introduce the Poisson exponential
\[ e(\fii;\cdot):=\sum_{n=0}^\infty \frac 1{n!}\, Q_n(\fii^{\ot
  n};\cdot)=I(\Exp\fii),\]
where
$$\Exp\fii=\bigg(\frac1{n!}\,\fii^{\ot n}\bigg)_{n=0}^\infty.$$
Then, one can  show that, for $\fii>-1$,
\begin{equation}\label{zebra}
e(\fii;\om)=\exp\big[\la\log(1+\fii),\om\ra-\la\fii\ra_{\tsigma}\big],\qquad
\om\in\OM. \end{equation}

\subsection{Second quantization on the marked Poisson space}

Let $B$ be a contraction on $L^2(\tsigma)$, i.e., $B\in {\cal
  L}(L^2(\tsigma),L^2(\tsigma))$, $\|B\|\le1$. Then, we can define the operator
  $\Exp B$ as the contraction on ${\cal F}(L^2(\tsigma))$ given by
\begin{align*}
&\Exp B \restriction {\cal F}_n(L^2(\tsigma)):=B\ot\dotsm\ot B\quad \text{($n $
  times)},\ n\in\N,\\
&\Exp B\restriction {\cal F}_0(L^2(\tsigma)):={\bf 1}.\end{align*}

For any selfadjoint positive operator $A$ in $L^2(\tsigma)$, we
have a contraction semigroup $e^{-tA}$, $t\ge0$, and it is
possible to introduce a positive selfadjoint operator $d \Exp A$
as the generator of the semigroup $\Exp(e^{-tA})$, $t\ge0$:
\begin{equation}\label{kon}
\Exp(e^{-tA})=\exp(-td\Exp A).\end{equation}
The operator $d\Exp A$ is called the second quantization of $A$. We denote by
$H^{\mathrm MP}_A$ the image of the operator $d\Exp A$ in the marked Poisson
space $L^2(\pi_\tsigma)$.

\begin{th}\label{th5.1}
Let ${\frak D}_0\subset\operatorname{Dom}A$\rom. Then\rom, the symmetric bilinear
form corresponding to the operator $H_A^{\mathrm MP}$ has the following
representation\rom:
\begin{equation}\label{5.18}
(H_A^{\mathrm MP}F_1,F_2)_{L^2(\pi_\tsigma)}=\int_{\OM}(\nMP F_1,A\nMP
F_2)_{L^2(\tsigma)}\,\pi_\tsigma(d\om) \end{equation}
for all $F_1,F_2\in\FPO.$\end{th}

\begin{rem}\label{rem5.1}
\rom{The bilinear form \eqref{5.18} uses the marked Poisson gradient $\nMP$ and
  a coefficient operator $A>0$. We will call
\[ {\cal E}^{\mathrm MP}_{\pi_\tsigma,A}(F_1,F_2)=\int_{\OM}(\nMP F,A\nMP
  G)_{L^2(\tsigma)}\,\pi _{\tsigma}(d\om)\]
the marked Poisson pre-Dirichlet form with coefficient $A$.}\end{rem}

\noindent {\it Proof of Theorem\/} 5.1. The proof is analogous to that of
Theorem~5.1 in \cite{AKR}. Using again the fact that ${\frak D}_0$ is an algebra
under pointwise multiplication, one easily concludes that, for any $F\in\FPO$ and
any $\om\in\OM$, the gradient $\nMP F(\om,(x,m))$ is
a function in ${\frak D}_0$ and hence
\[ (\nMP F,A\nMP G)_{L^2(\tsigma)}\in\FPO,\]
so that the form \eqref{5.18} is well-defined. Then, one verifies the formula
\eqref{5.18} by using Propositions~5.1, 5.2 and the explicit formula
for $d\Exp A$ on ${\cal F}_n({\frak D}_0)$:
\[ d\Exp A\,\fii^{\ot n}=n(A\fii)\hot \fii^{\ot(n-1)},\qquad \fii\in{\frak
  D}_0.\quad
\blacksquare\]

\subsection{The intrinsic Dirichlet operator as a second quantization}

The following two theorems are again analogous to ´the
corresponding results (Theorems~5.2 and~5.3) in \cite{AKR}, so we
omit their proofs.

Let us consider the special case of the second quantization operator $d\Exp A$
where  the operator $A$ coincides with the Dirichlet operator $H^{\XM}_\tsigma$.

\begin{th}\label{th5.2} We have the equality
\[ H^{\mathrm MP}_{\HH}=\HHH\]
on the dense domain $\FCpO$\rom. In particular\rom, for all $F_1,F_2\in\FCpO$
\begin{align*}
& \int_{\OM}\la \nabla^\Om F_1(\om),\nabla^\Om
F_2(\om)\ra_{T_\om(\OM)}\,\pi_\tsigma(d\om)\\
&\qquad = \int_{\OM}(\nMP F_1(\om),\HH\nMP
F_2(\om))_{L^2(\tsigma)}\,\pi_\tsigma(d\om),\end{align*}
or
\[ \nabla^{\Om*}\nabla^\Om=\nabla^{{\mathrm MP}*}\HH\nMP\]
as an equality on $\FCpO.$\end{th}

\begin{th}\label{th5.3}Suppose that the operator $\HH$ is essentially
  selfadjoint on the domain ${\frak
  D}_0\subset\operatorname{Dom}(H_\tsigma^\XM)$\rom. Then\rom, the intrinsic
  Dirichlet operator $\HHH$ is essentially selfadjoint on the domain
  $\FCO.$\end{th}

\begin{rem}\label{rem5.2}\rom{
Notice that in Theorem~\ref{th5.3} we do not suppose the operator $\HH$ to be
conservative. So, this theorem is a generalization of Theorem~\ref{th4.2} in the
special case where $\munu=\pi_\tsigma$.}\end{rem}

\begin{cor}\label{cor5.1}
Suppose that the condition of Theorem~\rom{\ref{th5.3}}
is satisfied and let $T^\Om_{\pi_\tsigma}(t)=\exp(-t\HHH)$\rom,
$t>0$\rom. Then\rom, for each $\fii\in{\frak D}_0$\rom, $\fii>-1$\rom, we have
\begin{equation}
\TTT\exp(\la\log(1+\fii),\cdot\ra)
=
\exp\big[\la\log(1+e^{-t\HH}\fii),\cdot\ra-\la(e^{-t\HH}-{\bf 1
})\fii\ra_\tsigma\big]. \label{mein}\end{equation}\end{cor}

\noindent {\it Proof}. The formula \eqref{mein} follows from
Proposition~\ref{prop5.1}, \eqref{zebra}, \eqref{kon} and Theorems~\ref{th5.2}
and~\ref{th5.3}.\quad$\blacksquare$

\begin{rem}\rom{
If $\HH$ is conservative, then
\[ \int(e^{-t\HH}-{\bf 1})\fii\,\,d\tsigma=0\qquad\text{for all }t\ge0,\]
and so in this case \eqref{mein} coincides with \eqref{peter}
for $\fii\in{\frak D}_0$, $\fii>-1$.}\end{rem}

\begin{center}\bf ACKNOWLEDGMENTS\end{center}

The authors were partially supported by the SFB 256, Bonn
University. Support by the DFG through Projects~436~113/39 and
436~113/43, and by the BMBF through Project~UKR-004-99 is
gratefully
 acknowledged.

\address{Institut f\"{u}r Angewandte Mathematik,
Universit\"{a}t Bonn, Wegelerstr. 6, D 53115 Bonn; and\\ SFB 256,
Univ.~Bonn; and\\
 CERFIM (Locarno); Acc. Arch. (USI); and\\
BiBoS, Univ.\ Bielefeld}\vspace{2mm}

\address{Institut f\"{u}r Angewandte Mathematik,
Universit\"{a}t Bonn, Wegelerstr. 6, D 53115 Bonn; and\\ SFB 256,
Univ.~Bonn; and\\ Institute of Mathematics, Kiev; and\\ BiBoS,
Univ.\ Bielefeld}\vspace{2mm}

\address{Institut f\"{u}r Angewandte Mathematik,
Universit\"{a}t Bonn, Wegelerstr. 6, D 53115 Bonn; and\\ BiBoS,
Univ.\ Bielefeld}\vspace{2mm}

\address{Radolfzellerstr.~9, D-81243 M\"unchen, Germany}

\end{document}